\input amstex.tex
\input amsppt.sty   
\magnification 1200
\vsize = 8.28 true in
\hsize=6.2 true in
\nologo
\NoRunningHeads        
\parskip=\medskipamount
        \lineskip=2pt\baselineskip=18pt\lineskiplimit=0pt
       
        \TagsOnRight
        \NoBlackBoxes

        \topmatter
         \title
      Infinite Energy quasi-periodic Solutions 
       \\to Nonlinear Schr\"odinger Equations on $\Bbb R$
        \endtitle
        
\author
         W.-M.~Wang \endauthor     
\address
{CNRS and D\'epartement de Math\'ematique, Universit\'e Cergy-Pontoise, 95302 Cergy-Pontoise Cedex, France}
\endaddress
        \email
{wei-min.wang\@math.cnrs.fr}
\endemail
\abstract
We present a set of {\it smooth infinite energy} global solutions (without spatial symmetry)
to the non-integrable, nonlinear Schr\"odinger equations on $\Bbb R$.
These solutions are space-time quasi-periodic with two frequencies each.
Previous results \cite{B2,~1}, and the generalization \cite {W2}, are quasi-periodic in time, but periodic in space.
This paper generalizes Bourgain's semi-algebraic set method \cite{B3} to analyze nonlinear PDEs,
in the {\it non-compact} space quasi-periodic setting on $\Bbb R$.
\endabstract
\keywords Quasi-momenta, Newton scheme, Anderson localization, small-divisors, rational approximation of diffeo-morphism, semi-algebraic geometry
\endkeywords
     
        \bigskip\bigskip
        \bigskip
        \toc
        \bigskip
        \bigskip 
        \widestnumber\head {Table of Contents}
        \head 1. Introduction to the Theorem
        \endhead
        \head 2. Green's function estimates in $(\theta, \phi)$
        \endhead
        \head 3. Nonlinear analysis -- proof of the Theorem 
        \endhead
        \endtoc
        \endtopmatter
        \vfill\eject
        \bigskip
       
\document
\head{\bf 1. Introduction to the Theorem}\endhead
We consider the nonlinear Schr\"odinger equation (NLS)  on $\Bbb R$:
$$
i\frac{\partial}{\partial t}u =-\frac{\partial ^2}{\partial x^2} u+|u|^{2p}u,\tag 1.1
$$
where $p\geq 1$ and $p\in\Bbb N$ is {\it arbitrary}; $u$ is a complex valued function on $\Bbb R\times \Bbb R$.
The purpose of this paper is to present possibly a new phenomenon, namely equation (1.1) admits families of global solutions, which
do not vanish at infinity, and have no apparent symmetry. 
Toward that end, we seek space-time quasi-periodic solutions of two frequencies each, i.e., $(2, 2)$-frequencies, in the form of a convergent sum:
$$
u(t, x)=\sum_{(n,j)\in\Bbb Z^2\times \Bbb Z^2}\frak a(n, j)e^{i(n\cdot\omega+M) t}e^{i(j\cdot \lambda+m) x}, \qquad \tag A
$$
where $\lambda=(\lambda_1, \lambda_2)\in (0, 2\pi]^2$ and $(m, M)\in (0, 2\pi]^2$ are the parameters; the coefficients $\frak a(n, j)\in\Bbb R$  and the frequency in time $\omega=(\omega_1, \omega_2)\in\Bbb R^2$ are to be determined for {\it appropriate}
$(\lambda, m, M)\in (0, 2\pi]^4\subset\Bbb R^4$. 

We view the series in (A) as a {\it quasi-periodic} Fourier series, and  $\frak a(n, j)$ the Fourier coefficients. 
In this terminology, the usual Fourier series is then a {\it periodic} 
Fourier series. We note, however, that the functions defined by  (A) {\it do not} decay in $x$ and a fortiori have {\it infinite} energy and mass for all time. 
Moreover, contrary to the space-periodic case, 
there is no reduced quotient space, namely the torus, on which $u$ has finite energy and mass.
So the problem considered here is a generalization of that in \cite{W2} on the torus. 

Our main result is 
that on a large (Cantor) set in $(\lambda, m, M)$, there are solutions of the form (A).
These solutions are bifurcations of solutions 
with {\it two} frequencies in time to the linear Schr\"odinger 
equation on $\Bbb R$, 
$$
i\frac{\partial}{\partial t}u =-\frac{\partial ^2}{\partial x^2} u;
$$
the $\omega$ in (A) are modulations of these two frequencies.

Below is the precise statement.
\bigskip
\proclaim{Theorem} 
Let $h_1$, $h_2\in\Bbb Z^2\backslash\{0\}$, $$h_1 \not\parallel h_2,$$ 
and $$u^{(0)}(t, x)=\sum_{k=1}^2 a_k e^{-i(h_k\cdot\lambda+m)^{2}t}e^{i(h_k\cdot \lambda+m)x}\tag U$$ 
be a solution to the linear Schr\"odinger equation, with
$\lambda=(\lambda_1, \lambda_2)\in (0, 2\pi]^2$, $m\in (0, 2\pi]$ and $a=(a_1, a_2)\in (0,\delta)^2$.   
For all $\epsilon>0$, 
there exists $\delta_0>0$, such that 
for all $\delta\in (0, \delta_0)$ and all fixed $a\in (0,\delta)^2$, there is a Cantor set $\Cal G\subset (0, 2\pi)^4$ with 
$$\text{meas }\Cal G/(2\pi)^4\geq 1-\epsilon,\tag 1.2$$ 
and a diffeomorphism: $(\lambda, m, M)\mapsto(\lambda, \omega)$ on $(0, 2\pi)^4$,
satisfying 
$$|\omega_k-(h_k\cdot\lambda+m)^2-M|=\Cal O(\delta^{2p}), \, k=1, 2.\tag 1.3$$
For all $(\lambda, m, M)\in\Cal G$, there is a space-time quasi-periodic solution to the NLS (1.1)
in the form (A), the Fourier coefficients $\frak a(n, j)$ decay sub-exponentially as $|(n,j)|\to\infty$, and 
$$|u(t, x)-\sum_{k=1}^2 a_k e^{i(-\omega_k+M)t}e^{i(h_k\cdot \lambda+m)x}|=\Cal O(\delta^p),\tag 1.4$$
for all $(t, x)\in\Bbb R\times \Bbb R$.
\endproclaim

\noindent{\it Remark.} Note that the Fourier phase $M$ does not appear in (U). This cancellation is due to the {\it first order} operator $i \partial/\partial t$. 
An equivalent formulation would be seeking solutions of the form 
$$
u(t, x)=\sum_{(n,j)\in\Bbb Z^2\times \Bbb Z^2}\frak a(n, j)e^{in\cdot\omega t}e^{i(j\cdot \lambda) x}, \qquad 
$$
to the {\it parameter dependent} equation: 
$$
(i\frac{\partial}{\partial t}+M)u =(i\frac{\partial}{\partial x}+m)^2u+|u|^{2p}u.
$$

\proclaim{Corollary} Assume that the hypothesis in the Theorem holds. There exists a set $\frak M\subset (0, 2\pi)^2$, satisfying 
$$\text{meas }\frak M/(2\pi)^2\geq 1-\epsilon,$$
such that for all $(m, M)\in\frak M$ 
the conclusion of the Theorem holds on a set in $\lambda$ (depending on $(m, M)$) of measure at least $(2\pi)^2(1-\epsilon)$; likewise, 
there exists a set $\Lambda\subset (0, 2\pi)^2$, satisfying 
$$\text{meas }\Lambda/(2\pi)^2\geq 1-\epsilon,$$
such that for all $\lambda\in\Lambda$ 
the conclusion of the Theorem holds on a set in $(m, M)$ (depending on $\lambda$) of measure at least $(2\pi)^2(1-\epsilon)$.
\endproclaim

\demo{Proof} This follows from the Theorem by using the Fubini Theorem.
\hfill $\square$
\enddemo
\bigskip
\subheading{1.1 Some background}

It is well known that the Cauchy problem for the NLS in (1.1), is globally well-posed  
in $H^1(\Bbb R)$, using Sobolev embeddings and energy conservation. (For $1\leq p\leq 2$, (1.1) is, in fact, globally well-posed
in $L^2(\Bbb R)$ \cite{Tsu}\cite{Do}.) For periodic $u$, set $$\Bbb T:=2\pi\Bbb R/\Bbb Z,$$
then (1.1) is globally well-posed in $H^1(\Bbb T)$. Periodic $u$  are not localized in space, but 
on the quotient space there is finite $H^1$ norm. 

It is natural to seek quasi-periodic in space solutions $u$ to (1.1).
If they existed, such $u$, a fortiori, would have both infinite energy and infinite mass ($H^1$ and $L^2$ norms). 
Since there is no lack of regularity, these $u$ present a different phenomenon from that of rough solutions, $u\in H^s(\Bbb R)$ for $s<1$. 
There are, however, few such known results.

In \cite{DG}, Damanik and Goldstein constructed analytic small amplitude space quasi-periodic and time almost-periodic Cauchy solutions to KdV equations. 
Such solutions were previously constructed in \cite{BdME} for the cubic NLS, i.e.,  (1.1) when $p=1$, under certain spectral assumptions. Both results seem to hinge on integrability. 
The paper \cite{DG} uses inverse spectral analysis, and detailed spectral information of the underlying linear quasi-periodic Schr\"odinger operator in one dimension.
In particular, it requires that the quasi-periodic Schr\"odinger operator has purely absolutely continuous spectrum.  
Point spectrum is excluded. In the case of the non-integrable NLS in (1.1), however, we cannot exclude point spectrum,
for the linearized space-time quasi-periodic operators, which are in two dimensions (space plus time).
The proof of the Theorem does not rely on spectral method. 

In the non-integrable case, Oh proved local well-poseness of Cauchy problems for almost-periodic initial data $u_0$, and global well-posedness 
for limit-periodic $u_0$, for (1.1), in \cite{O1, 2}. (Limit periodic functions are functions with {\it one} basic period and all other periods
are rational multiples of it. In particular, it cannot be quasi-periodic, which has at least {\it two} basic periods.) Dodson, Soffer and Spencer \cite{DoSS}
proved local well-posedness for real analytic data, and global existence for a regularized NLS. For NLS on the lattice $\Bbb Z$, they proved that 
solutions are bounded by a fractional power in time.  
We also mention that \cite{CHKP} has a result on finite energy perturbation of a periodic NLS. 

\bigskip
\subheading{1.2 About the Theorem}

In the Theorem and the Corollary, we prove existence of space-time quasi-periodic solutions with two frequencies each, to (1.1) for all $p\geq 1$. 
These solutions appear to be first instances of global in time, smooth, non-localized solutions (and without symmetry), to non-integrable equations such as (1.1). It is a set of solutions
parametrized by four parameters, namely the two space frequencies $\lambda=(\lambda_1, \lambda_2)$, and the two Fourier phases, the quasi-momenta, $m$ and $M$.   

The method devised in this paper is general and applicable in any dimensions, and any number of frequencies in time,
as long as they are {\it independent} from the space frequencies $\lambda$, and are of order $1$. (Cf. also \cite{W3}.) 
For example, it is applicable to {\it forced} equations. The forcing frequencies would then provide the aforementioned additional $\Cal O(1)$ parameters.  

For autonomous equations like (1.1), however,  aside from $\lambda$, the {\it only} other natural $\Cal O(1)$ parameters in the problem are, the phases $m$ and $M$. The result in the paper 
therefore pertains to two basic frequencies in time. It is in contrast to the space periodic case in \cite{W2},
where time quasi-periodic solutions with arbitrary number of frequencies are constructed, by using 
amplitude-frequency modulations, which yield $o(1)$ parameters. 
\bigskip
Let us elaborate the reasons for the difference.
\smallskip
\noindent$\bullet$ {\it The periodic Laplacian}
\smallskip
The spectrum of the Laplacian on the torus $\Bbb T^d$ is the set:
$$\sigma(-\Delta)=\{\sum_{i=1}^d j_i^2:=|j|^2  |\, j_i\in\Bbb Z, i =1, 2, ...,d\}\subseteq\Bbb Z^+.$$
The operator is {\it elliptic}: as $|j|\to \infty$,  $|j|^2\to \infty$, and modulo degeneracies, the spectral gaps are at least $1$.
\bigskip
\noindent$\bullet$ {\it The quasi-periodic Laplacian}
\smallskip
The Laplacian here, on the other hand, acts on space quasi-periodic Fourier series.
Using the Ansatz (A), it gives instead
$(j\cdot\lambda +m)^2$, which could be dense on $\Bbb R^+$,
and arbitrarily close to $0$, as $|j|\to\infty$. So it
is {\it non-elliptic}, and there is {\it no} spectral gap. 

One may also amplitude-frequency modulate in the present setting. However, since the 
nonlinear term is non-decaying and there is no spectral gap (dense spectrum), the extracted 
$o(1)$ parameters are too small to control the small divisors; while in \cite{W2} they suffice. This is because 
once the degeneracies are taken into account, the spectral gaps compensate for the non-decaying nonlinear term.

\noindent{\it Remark.} In \cite{BW}, the spectrum (of the random Schr\"odinger operator) is also dense, but
the exponentially decaying nonlinear term makes up for this deficiency.
\smallskip
Returning to the problem at hand, like the space direction, the time direction is also non-elliptic, with $i\partial/\partial t$ yielding $(n\cdot\omega+M)$. This leads to new types of small divisors,
which we analyze in the semi-algebraic setting. In fact, semi-algebraic
geometry is brought in precisely to control non-ellipticity in dimensions two and above, as there remains a notion 
of {\it steepness}, when $|(n, j)|\gg 1$.

In general, {\it parameters} are needed for this type of constructions.
Otherwise the resonances could possibly conspire to
growth of Sobolev norms. See, for example,  \cite {CKSTT} \cite{GHP} and references therein
for finite time growth examples on the torus.
In this paper, {\it phase-frequency modulation} and semi-algebraic geometry permit, to select the parameters 
to control the dense spectrum, and prove the Theorem.

Without such selections, 
one would not have expected
space-time quasi-periodic solutions, for non-integrable equations such as (1.1). 
Indeed, as mentioned earlier, there is a dearth of global, non-localized solutions without spatial symmetry, to 
non-integrable, nonlinear PDEs.
\smallskip
\noindent {\it Remark.} For the nonlinear wave equations on $\Bbb R$,  using finite speed of propagation, one may show 
that for bounded smooth initial data, the solutions are bounded by a fractional power in time \cite{DoSS}. 
Due to infinite speed of propagation, however, there are no such known apriori bounds for 
the nonlinear Schr\"odinger equations on $\Bbb R$.
\bigskip
\subheading{1.3 Anderson localization and semi-algebraic set method}
\smallskip
\noindent$\bullet$ {\it The linear problem}
\smallskip
We use the Anderson localization approach to construct the space-time quasi-periodic solutions in the Theorem and Corollary. 
Here, Anderson localization is in the {\it Fourier space}. 
The mathematical theory of Anderson localization debuted in the seminal paper of Fr\"ohlich and Spencer \cite{FS}. In the paper, they show how
to probabilistically control small divisors and prove eigenfunction localization of random linear systems (random Schr\"odinger operators). 
This novel mechanism is very flexible and applicable to quasi-periodic systems as well, see the application in one dimension \cite{FSW}.
\smallskip
Quasi-periodic Anderson localization in higher dimensions, is however, considerably more difficult.  This is because, in one dimension,
the dynamical variable, seen as the ``time", is a {\it scalar}; while in higher dimensions, it is a {\it vector}. To deal with this new situation, 
Bourgain introduced novel techniques from harmonic analysis and semi-algebraic geometry. In the breakthrough paper \cite{BGS}, Bourgain, Goldstein and
Schlag proved, for the first time, Anderson localization for quasi-periodic systems on the $\Bbb Z^2$ lattice. Later in \cite{B3}, Bourgain
generalized the result to $\Bbb Z^d$, for arbitrary $d$, by developing powerful techniques of 
semi-algebraic variable reductions. (Cf. also the dissected, recent rendition \cite{JLS}.)

Finally in \cite{BK}, Bourgain and Kachkovskiy treated a certain degenerate quasi-periodic system of two identically equal frequencies,
on the $\Bbb Z^2$ lattice. (Here by degenerate, we mean that the number 
of free parameters, the {\it independent} quasi-periodic frequencies, is strictly less than the dimension of the lattice.) 
\bigskip

\noindent$\bullet$ {\it The nonlinear problem}
\smallskip
Parallel to the linear Anderson localization, there is the {\it nonlinear} Anderson localization. 
This theory was developed by Bourgain to construct time quasi-periodic, and {\it space-periodic} solutions to 
{\it parameter-dependent} nonlinear PDEs. The application of Anderson localization
to construct time periodic solutions was first made by Craig and Wayne \cite{CW}. 
It was recently generalized to the space quasi-periodic setting in \cite{W3}. The bridge between
the nonlinear and the linear theory is provided by a Newton scheme -- after each linearization, 
one falls back on the linear analysis. 

Chaps.~18-19 in \cite{B2} recount this nonlinear theory, see also the generalizations in \cite{W2}, and the review \cite{W4}. 
We note that \cite{W2} treats the {\it original} nonlinear equations, such as those in (1.1); while \cite{B1,~2}, only 
a family of parameter-dependent equations. (As mentioned earlier, \cite{W2} deals with the original equations by {\it extracting} parameters from
the nonlinear term. The extraction method is algebraic.) 
This theory was intended for Hamiltonian PDEs 
on {\it compact} manifolds, mainly represented by flat tori. The PDEs can then be seen as infinite dimensional 
systems of ODEs. Cf.  \cite{EK}\cite{PP} for related KAM-type theory.

The nonlinear Anderson localization theory is, however, more general. 
It is not limited to the compact setting. In particular, it is {\it not} restricted to the space-periodic case, 
namely the tori. One can {\it equally well} pose the problem in the {\it non-compact} space quasi-periodic setting on $\Bbb R$, by using instead
space-time quasi-periodic Fourier series. 
The Newton scheme then enables us to use the arsenals from the
linear Anderson localization theory for quasi-periodic systems.  

\noindent{\it Remark.} Unlike the space-periodic nonlinear Schr\"odinger equations, the space quasi-periodic 
nonlinear Schr\"odinger equations do not seem to have a natural dynamical systems interpretation. 

However, since the space quasi-periodic setting is non-compact, the Laplacian is non-elliptic, as mentioned earlier. 
The resulting quasi-periodic system is at least in two dimensions. 
So we use semi-algebraic geometry, so far, the only available method to deal with quasi-periodic systems in dimensions two and above.

Previously, in the space-periodic case in Chaps.~19 \cite{B2} and \cite{W2}, 
the linear analysis is done in the {\it flat} parameter space. Here it would be $(\lambda, \omega, m, M)\in\Bbb R^6$. 
Because of the non-elliptic Laplacian, the linear analysis here, which partly resorts to semi-algebraic geometry, is considerably more delicate.
The dimension of the (independent) parameter set, here four, plays an essential role. Consequently, the analysis
restricts to the neighborhood of the co-dimension $2$ sub-manifold with coordinates $(\lambda, \omega (\lambda, m, M), m, M)$.

One of the key new aspect is, to make rational approximations of diffeomorphisms, on open (sub)sets in $(\lambda, m, M)$. The union of these open sets eventually turns into the 
Cantor set $\Cal G$ in the Theorem. The other new aspect, is that since {\it both} space and time 
are non-compact, we need to generalize the semi-algebraic analysis in \cite{B3}, which is in the compact 
setting, by developing a (weakly) nonlinear Lipschitz eigenvalue analysis.


In the nonlinear analysis, aside from the measure, the {\it structure} of the subset in $(\lambda, m, M)$, on which one may do linear
analysis, comes into play. The fact that this subset maybe taken to be a set of open sets of certain size, together with the double exponential convergence 
of the Newton scheme, enable space-time quasi-periodic solutions. 
\smallskip
In \cite{W3}, we constructed space quasi-periodic standing waves on $\Bbb R^d$, for arbitrary $d$, using semi-algebraic geometry. In that case,
the nonlinearity $|u|^{2p}$ is independent of time, and the quasi-periodic frequencies $\lambda$ are the given parameters. 
One works in the given flat parameter space. There is no need to change variables.

In this paper, the space frequencies $\lambda=(\lambda_1, \lambda_2)$, and the phases $m$ and $M$, are the given parameters;
the time frequencies $\omega=(\omega_1, \omega_2)=\omega (\lambda, m, M)$ are, however, only {\it functions} of the given parameters $\lambda, m$ and $M$.
The parameters $(\lambda, \omega)$ describe the {\it frequency manifolds}. As mentioned earlier, these manifolds need to be close to algebraic ones. 
In the $(\lambda, \omega)$ variables, the analysis is related to that in \cite{W3}.  
\smallskip
To conclude, using semi-algebraic geometry, this paper provides a possibly new way  
to deal with non-decaying nonlinearity, in the non-compact and non-integrable setting.  
\bigskip

\noindent{\it Remark.} In \cite{GYZ}, time quasi-periodic solutions were constructed for a quasi-periodic nonlinear 
lattice Schr\"odinger equation using a KAM method, see also \cite{Y}. However, since the authors perturb about localized solutions, the diagonal nonlinear term $|u_j|^{2}u_j$
decays rapidly on $\Bbb Z$, and the solutions $u$ have {\it finite} energy and mass. The result is instead, in the spirit of \cite{BW}, which uses nonlinear Anderson localization.
\bigskip
\noindent{\it Organization of the paper.}

In sect.~2, we use the Ansatz in (A), and make a Lyapunov-Schmidt decomposition of the resulting equations 
into the $P$ and $Q$-equations. The $P$-equations are used to solve for the
amplitudes $\frak a(n, j)$; while the $Q$-equations, the modulated frequencies $\omega_1$ and $\omega_2$, keeping $a_1$ and
$a_2$ fixed. The $P$-equations are solved using a Newton iteration. Sect.~2
analyzes the ensuing linearized operators. Sect.~3 uses the linear 
estimates, to solve the $P$ and $Q$-equations and prove the Theorem. 
\bigskip
\noindent{\bf Acknowledgement.} The Theorem answers a question raised by T. Spencer; the author
would like to thank him for discussions and comments. She also thanks Wencai Liu for proofreading and suggestions; and NYU Shanghai for hospitality in 
spring 2019.

\head {\bf 2. Green's function estimates in $(\theta, \phi)$}\endhead
\smallskip
\subheading {2.1 The Lyapunov-Schmidt decomposition}
\smallskip
Using the ansatz in (A) to express (1.1); by analogy with the standard Fourier series, write $\hat u$ for $\frak a$ and define $\hat {\bar u}$ to be
$\hat {\bar u}(n, j)=\bar \frak a(-n, -j)$ for all $(n,j)\in\Bbb Z^{2}\times\Bbb Z^2$. To simplify notations, write $\hat v$ for $\hat {\bar u}$.
Equation (1.1) can then be written as a nonlinear (infinite) matrix equation: 
$$\text{diag }(n\cdot\omega+M+(j\cdot\lambda+m)^2)\hat u+(\hat u*\hat v)^{*p}* \hat u=0,\tag 2.1$$
where diag $\cdot$  denotes a diagonal matrix, $(\lambda=(\lambda_1, \lambda_2), m, M)\in (0, 2\pi]^4\subset\Bbb R^4$
are the free parameters, and $\omega=\omega (\lambda, m, M)\in\Bbb R^2$ is to be determined. 

We work with (2.1), for simplicity, drop the hat and write $u$ for $\hat u$ and $v$ for $\hat v$ etc.
We seek solutions $u=u(\lambda, m, M)$ (and $v=v(\lambda, m, M)$) close to the $(2, 2)$-frequencies linear solution $u^{(0)}$ in (U) with time-frequencies: 
$$\omega^{(0)}=\{(h_1\cdot\lambda+m)^2+M, (h_2\cdot\lambda+m)^2+M\},\, (h_1\not\parallel h_2\neq 0),\tag 2.2$$
and space-frequencies: $\{h_1\cdot\lambda, h_2\cdot\lambda\}$, 
and small amplitudes $a=\{a_1, a_2\}$ satisfying $\Vert a\Vert~=\Cal O(\delta)\ll 1$.
Comparing (U) with (A), one may write the Fourier support of $u^{(0)}$ as
$$\text{supp } {u}^{(0)}=\{(-e_{1}, h_1), (-e_{2}, h_2)\},\tag 2.3$$ 
where $e_1$ and $e_2$ are the two base vectors in $\Bbb Z^2$ in the $n$-direction.

It is convenient to rescale: 
$$\aligned &u\to \delta u\\
& a_k \to \delta a_k,\, k=1,2. \endaligned \tag R$$
So in what follows, the ``new" $a_k=\Cal O(1)$, $k=1, 2$, and since they are fixed, we shall not make explicit the dependence on $a$.
(In the Theorem and Corollary, however, we revert to the ``original" $a_k, \, k=1,2$.)
Completing (2.1) by the equation for the complex conjugate, we arrive at the system on the lattice $$\Bbb Z^{2}\times\Bbb Z^2\times\{0, 1\}:=\Bbb Z^4\times\{0, 1\}:$$
$$
\cases
\text{diag }(n\cdot\omega+M+(j\cdot\lambda+m)^2)u+\delta^{2p}(u*v)^{*p}* u=0,\\
\text{diag }(-n\cdot\omega+M+(j\cdot\lambda+m)^2)v+\delta^{2p}(u*v)^{*p}* v=0.
\endcases\tag 2.4
$$
By ``supp", we mean the Fourier support, so write $\text{supp } u^{(0)}$ for 
$\text{supp } {\hat u}^{(0)}$ etc. Let 
$$\aligned\Cal S=&\text{supp } u^{(0)}\oplus\text{supp } {\bar u}^{(0)}\\
:=&S\oplus\bar S.\endaligned\tag 2.5$$

Denote the left side of (2.4) by $F(u, v)$. We make a Lyapunov-Schmidt decomposition into the $P$-equations:
$$ F(u, v)|_{\Bbb Z^{4}\times\{0, 1\}\backslash\Cal S}=F(u, v)|_{\Cal S^c}=0,\tag 2.6$$
and the $Q$-equations:
$$ F(u, v)|_{\Cal S}=0,\tag 2.7$$
which are solved separately. We seek solutions such that 
$$(u, v)|_{\Cal S}=(u^{(0)}, v^{(0)}).$$ 

The $P$-equations are infinite dimensional and determine $(u, v)$ in $\Cal S^c$; 
the $Q$-equations are $4$ dimensional and determine the frequency $\omega=\{\omega_k\}_{k=1}^2$. (Due to 
symmetry, two of the four $Q$-equations are identical.)
\bigskip
\noindent $\bullet$ {\it The $P$-equations}
\smallskip
We use a Newton scheme to solve the $P$-equations. 
Let $F'$ be the linearized operator
on $\ell^2(\Bbb Z^{4})\times\ell^2(\Bbb Z^{4})$,
$$F'=D+\delta^{2p}H,\tag 2.8$$
where 
$$
D =\pmatrix \text {diag }(n\cdot\omega+M+(j\cdot\lambda+m)^2)&0\\ 0& \text {diag }(-n\cdot\omega+M+(j\cdot\lambda+m)^2)\endpmatrix,\tag 2.9$$
and
$$\aligned
H&=\pmatrix (p+1)(u*v)^{*p}*& p(u*v)^{*p-1}*u*u*\\ p(u*v)^{*p-1}*v*v*& (p+1)(u*v)^{*p}*\endpmatrix
\quad  (p\geq 1),\\
\endaligned\tag 2.10$$
is a convolution matrix. $F'$ is to be evaluated near $\omega=\omega^{(0)}$, $u=u^{(0)}$, $v=v^{(0)}$.

Recall next the formal Newton scheme:
$$\Delta \pmatrix u\\v\endpmatrix=-[F'_{\Cal S^c}(u, v)]^{-1} F(u,v)|_{\Cal S^c},$$ 
where the left side denotes the correction to $\pmatrix u\\v\endpmatrix$, $F'_{\Cal S^c}(u, v)$ is the linearized operator evaluated at $(u, v)$: $F'(u, v)$, and restricted to $\Cal S^c$: 
$F'_{\Cal S^c}(u, v)(x, y)= F'(u, v)(x,y)$, for $x, y\in \Cal S^c$; likewise $F(u, v)|_{\Cal S^c}$ is $F(u, v)$ restricted to ${\Cal S^c}$:
$[F(u, v)|_{\Cal S^c}](x)=F(u, v)(x)$ for $x\in \Cal S^c$. 
\bigskip
Since we seek solutions close to $(u^{(0)}, v^{(0)})$, which has compact support in $\Bbb Z^4\times \{0, 1\}$, 
we adopt a {\it multiscale} Newton scheme as follows:

\noindent At iteration step $(i+1)$, choose an appropriate scale $N$ and estimate $[F'_N]^{-1}$, 
where $F'_N$ is $F'$ restricted to 
$$[-N, N]^{4}\times \{0, 1\}\backslash \Cal S\subset \Bbb Z^4\times \{0, 1\},$$
and evaluated at $u^{(i)}$ and $v^{(i)}$:  $F'_N=F'_N(u^{(i)}, v^{(i)})$.
We call $[F'_N]^{-1}$, the Green's functions. Define the $(i+1)$-th correction to be:  
$$\Delta \pmatrix u^{(i+1)}\\v^{(i+1)}\endpmatrix=-[F'_{N}(u^{(i)}, v^{(i)})]^{-1} F(u^{(i)},v^{(i)}),\tag 2.11$$ 
and 
$$\aligned u^{(i+1)}&=u^{(i)}+\Delta u^{(i+1)},\\
v^{(i+1)}&=v^{(i)}+\Delta v^{(i+1)},\endaligned\tag 2.12$$
for all $i=0$, $1$, $2$, ... 
\bigskip
\noindent $\bullet$ {\it The $Q$-equations}
\smallskip
The $Q$-equations in (2.7) determine the frequencies $\omega=\{ \omega_k\}_{k=1}^2$:
$$\omega_k=(h_k\cdot\lambda+m)^2+M+\delta^{2p}\frac{[{(u*v})^{*p}*u] (\lambda, m, M)}{a_k}(-e_k, h_k),\, k=1, 2.\tag 2.13$$
(Recall that $u, v$ also depend on $a$, but is not made 
explicit since $a$ is fixed.) The above equations, 
are solved exactly, so $F_\Cal S=0$ always.

\noindent{\it Remark.} Generally speaking, the $Q$-equations are solved using implicit function theorem, cf. \cite{B1,~2}\cite{W2}. However,
when the nonlinear analysis uses semi-algebraic geometry in sects.~3.5 and 3.6, we will need to be more precise, and make rational approximations
of this diffeomorphism beforehand, in sect.~3.4.
\bigskip
\subheading {2.2 Invertibility of the linearized operators -- the $(\theta, \phi)$ estimates}

We 
add a {\it two-dimensional} auxiliary variable $(\theta, \phi)\in\Bbb R^2$:
$$F'(\theta, \phi)=D(\theta, \phi) +\delta^{2p} H,\tag 2.14$$ 
where 
$$D(\theta, \phi) =\pmatrix D_+&0\\ 0& D_-\endpmatrix,\tag 2.15$$
and 
$$D_{\pm}:=\text {diag }[\pm (n\cdot\omega+\theta)+M+(j\cdot\lambda+\phi+m)^2],\tag 2.16$$
$H$ as in (2.10), and with a slight abuse, used the same symbols $F'$ and $D$.  We estimate
 the truncated Green's functions $[F'_N(\theta, \phi)]^{-1}$ in $(\theta, \phi)$ for all $N$,
and eventually use the covariance with respect to the $\Bbb Z^4$ action on $\Bbb R^2$: 
$$(\theta, \phi)\mapsto (\theta+n\cdot \omega, \phi+j\cdot\lambda),\tag 2.17$$ to deduce estimates for 
$[F'_N(\theta=0, \phi=0)]^{-1}:=[F'_N]^{-1}$, the Green's functions used in the Newton scheme (2.11)-(2.12).

In \cite{W2}\cite{B1,~2}, the space direction is compact (namely a torus), there is only the auxiliary variable $\theta$, which
can be bounded in terms of the (Fourier) scales $N$ to prove invertibility. (In \cite{W3}, the $\theta_i$ are also bounded by the scales $N$, as one
seeks nonlinear eigenvalues.) Here, however, since one works on $\Bbb R\times\Bbb R$, there is a new phenomenon, namely 
$(\theta, \phi)$ cannot be bounded in terms of the scales $N$.  The division of labor is that semi-algebraic geometry takes 
care of the region of $(\theta, \phi)$ bounded in terms of $N$, in sects.~2.5-2.8; while Lipschitz eigenvalue analysis
does the (unbounded) complement, in sect.~2.4. 

Denote the linearized operator $F'$ by $T$; and $F'_N$, $T_N$. The goal of this
section is 
to estimate the Green's functions $T_N^{-1}(\theta, \phi)$, 
away from a set in $(\theta, \phi)$ of small sectional measure. Since $u$ and $v$ will be determined iteratively in sect.~3, they depend on the 
scales $N$. We denote them by $u_N$ and $v_N$ in this section. 

We do not use the relation (2.13) here, and make separate assumptions on 
$u_N(v_N)$,  and $\lambda, m, M$ and $\omega$.
These assumptions will be verified in the nonlinear analysis in sect.~3.

\noindent{\it Remark.} It should be noted, however, that $\theta$ and $\phi$ are unbounded for the linear analysis {\it only}.
For the nonlinear analysis, when using the covaraince relation (2.17) and semi-algebraic projections,
$\theta$ and $\phi$ are always {\it bounded} by the scales in the Newton iteration.

\bigskip
We assume $0<\delta\ll 1$ sufficiently small, and begin with the initial scales.
\bigskip
\subheading {2.3  The $(\theta,\phi)$ estimates for initial scales}

Let $s>1$ and $\tau>0$, satisfying $0<s\tau<1$. Set the initial scales $N\in\Bbb N$ to be:
\smallskip
\noindent {\bf {Initial scales:}}

$$N\in [(\log|\log\delta|)^{1/\tau}, |\log\delta|^s].\tag 2.18$$
Set $$N_0=(\log|\log\delta|)^{1/\tau}.\tag 2.19$$

\noindent For the initial scales, we may assume $\lambda,\omega, m$ and $M$ are independent variables, and 
work in $(\lambda,\omega, m, M)\in \Sigma\subset\Bbb R^6$, where $\Sigma$ is a bounded set.

\noindent {\bf {Property V1:}}

Let $K>1$, $\tilde\gamma>\gamma>0$.
We say a real sequence:  $$w_N:\, \Sigma\mapsto \ell^2(\Bbb Z^{4}),\, N\in\Bbb N$$
has Property V1 at $(\lambda,\omega, m, M)\in \Sigma$ , if 
$$w_N(k)=0,\text{ for } k\notin[-N^K, N^K]^{4},\, K>1,$$
and
$$|w_N(k)|\leq e^{-|k|^\gamma},\, k\in\Bbb Z^{4} \, (0<\gamma<1),\tag \dag$$
$$\Vert w_N-w_{N+1}\Vert_{\ell^2(\Bbb Z^{4})} \leq  e^{-\tilde\gamma N}\,(\tilde\gamma>0),\tag\dag\dag$$
where $w_N(k)$ stands for $w_N(\lambda,\omega, m, M)(k)$.
\bigskip
\proclaim {Lemma A} Let $0<\delta\ll 1$, and consider the initial scales $N$ in (2.18).
For any $(\lambda, \omega, m, M)\in \Sigma$, if $u_N$ and $v_N$ have Property V1,  then there exists a subset $\Theta_N\subset \Bbb R^2$, whose sectional measures
satisfy
$$\text{meas } [\theta| \text {fixed }\phi, (\theta, \phi)\in \Theta_N] \leq e^{-N^\tau}\, (\tau>0),\tag 2.20$$
and $$\text{meas } [\phi| \text {fixed }\theta, (\theta, \phi)\in \Theta_N]\leq e^{-N^\tau} \, (\tau>0).\tag 2.21$$
If $(\theta,\phi)\notin \Theta_N$, the linearized operator $F_N':=T_N$ in (2.11), with
$u^{(i)}=u_N$ and $v^{(i)}=v_N$, satisfy the estimates
$$\Vert [T_N(\lambda,\omega, m, M; u_N, v_N)(\theta,\phi)]^{-1}\Vert_{\text{Op}}\leq e^{N^\sigma}\, (0<\tau<\sigma<1),\tag 2.22$$
$$|[T_N({\lambda, \omega, m, M; u_N, v_N})(\theta,\phi)]^{-1}(k,k')|\leq e^{-|k-k'|^\gamma}\, (0<\tau<\sigma<\gamma<1), \forall |k-k'|>N/10.\tag 2.23$$
\endproclaim

\demo{Proof}
We perturb about the diagonals: 
$$D_{\pm}:=\text {diag }[\pm (n\cdot\omega+\theta)+M+(j\cdot\lambda+\phi+m)^2].$$
It suffices that 
$$|D_{\pm}(\theta, \phi)|>\delta^{2p-1},$$
for $(n, j)\in[-N, N]^4$. Resolvent series in $\delta$, Property V1 and the relation (2.18), lead to (2.20)-(2.23), if $0<s\tau<1$ and $s\sigma>1$. 
Note that this holds for {\it all}  $(\lambda, \omega, m, M)\in \Sigma$, provided Property V1 is satisfied.
\hfill$\square$
\enddemo
\bigskip
\noindent  {\bf The regions (i) and (ii)}
\smallskip
Proceeding to larger scales $N$, $N\geq |\log\delta|^s$ $(s>1)$, there are two cases: 
\item {(i)} $|\theta|\leq 2e^{2(\log N)^2}$ and $|\phi|\leq 2e^{(\log N)^2}$;
\item{(ii)} otherwise.

\noindent (In fact sufficiently high degree polynomials in $N$ suffice for the above division.)
\smallskip

We start with the easier region (ii).
\bigskip
\subheading {2.4  The $(\theta,\phi)$ estimates in region (ii)}

We continue to view $\omega$ as an independent variable, and work in 
$$(\lambda, \omega, m, M)\in \Sigma\subset\Bbb R^6,$$
where $\Sigma$ is bounded and open.
For 
$$w:\, \Sigma\mapsto \ell^2(\Bbb Z^{4}),$$
with $w(\lambda, \omega, m, M)(k):=w(k)\neq 0$, for $k$ in a bounded subset of $\Bbb Z^4$, and 
$0$ otherwise,
define
$$\Vert  \text {Lip }w\Vert_{ \ell^2(\Bbb Z^{4})}:=\sup_x\Vert  \text {Lip}_x w\Vert_{ \ell^2(\Bbb Z^{4})},$$ 
where  $x$ stands for $\lambda, \omega, m$ or $M$.
\smallskip
\noindent {\bf {Property V2:}}

We say a real sequence $w_N$: 
$$w_N:\, \Sigma\mapsto \ell^2(\Bbb Z^{4}),$$
has Property V2, if it satisfies Property V1 on $\Sigma$ and 
$$\Vert \text {Lip }w_N\Vert_{\ell^2(\Bbb Z^{4})}=\Cal O(1),$$
uniform in $N$.

\bigskip
\proclaim{Lemma B} Let $0<\delta\ll 1$, $\lambda$ and $\omega$ be Diophantine, satisfying
$$\aligned \Vert j\cdot\lambda\Vert _{\Bbb T}&\geq \frac{\kappa}{|j|^{\rho}},\, 0<\kappa <1,\, \rho>3,\\
\Vert n\cdot\omega\Vert _{\Bbb T}&\geq  \frac{\kappa}{|n|^{\rho}}, \, 0<\kappa <1,\, \rho>3,\endaligned$$
for $j\neq 0$ and $n\neq 0$, where $\Vert \, \Vert _{\Bbb T}$ denotes distance to the integers. Assume $u_N$ and $v_N$ satisfy Property V2.
In region (ii), for any fixed $\phi$, there is a family of Lipschitz functions:
$\theta_i$, $i=1, 2, ...,  2(2N+1)^4$, satisfying
$$\Vert \theta_i\Vert _{\text{Lip}(\omega)}=\Cal O(N),\,  \Vert \theta_i\Vert _{\text{Lip}(M)}=\Cal O(1),$$
$$\Vert \theta_i\Vert _{\text{Lip}(\lambda)},  \Vert \theta_i\Vert _{\text{Lip}(m)}=\Cal O(N|\phi|),\tag 2.24$$
such that if
$$|\theta-\theta_i|>2e^{-N^{\sigma'}}, \,  0<\sigma'<\sigma,$$
then (2.22) and (2.23) are satisfied, for all $N\gg 1$. This excises a set in $\theta$ of measure satisfying (2.20). 
Likewise,  
for any fixed $\theta$, there is a family of Lipschitz functions $\phi_i$,  $i=1, 2, ...,  2(2N+1)^4$, satisfying 
$$\Vert \phi_i\Vert _{\text{Lip}(\omega)}=\Cal O(Ne^{-(\log N)^2}),\,  \Vert \phi_i\Vert _{\text{Lip}(M)}=\Cal O(e^{-(\log N)^2}),$$
$$\Vert \phi_i\Vert _{\text{Lip}(\lambda)}=\Cal O(Ne^{-(\log N)^2}\sqrt {|\theta|}),  \Vert \phi_i\Vert _{\text{Lip}(m)}=\Cal O(e^{-(\log N)^2}\sqrt {|\theta|}),\tag 2.25$$
such that if
$$|\phi-\phi_i|>2e^{-N^{\sigma'}},\, 0<\sigma'<\sigma,$$
then (2.22) and (2.23) are satisfied, for all $N\gg1$, leading to excision of a set in $\phi$ satisfying (2.21). 
\endproclaim
\smallskip 
\noindent{\it Remark.}  The lowering of $\sigma$ to $\sigma'$ is for the sake of (2.23). Note that in Lemma~A, one may take as Lipschitz functions the zeroes of the $D_\pm$, i.e., $\theta$ and $\phi$ 
such that $D_\pm=0$. 

\demo{Proof} There are three sub-cases:
\item{(a)} $|\theta|>2e^{2(\log N)^2}$, $|\phi|\leq e^{(\log N)^2}$; 
\item{(b)} $|\theta|>2e^{2(\log N)^2}$, $|\phi|>e^{(\log N)^2}$;
\item{(c)} $|\theta|\leq 2e^{2(\log N)^2}$, $|\phi|>2e^{(\log N)^2}$.
\bigskip 
\noindent Case (a). Write $T_\pm(n, j)$ for the diagonal elements of $T_N$: 
$$T_\pm(n, j)=D_\pm(n, j)+\delta^{2p}h=\pm(n\cdot\omega+\theta)+(j\cdot\lambda+\phi)^2+\delta^{2p}h,$$
where $h$ is the diagonal of $H$. For any fixed $\phi\in[-e^{(\log N)^2}, e^{(\log N)^2}]$, since $|\theta|>2e^{2(\log N)^2}$,
$$|T_\pm(n, j)|\geq \frac{1}{2}e^{2(\log N)^2},\tag *$$
for all $(n, j)\in [-N, N]^4$. 
So $$\Vert T_N^{-1}\Vert \leq 3 e^{-2(\log N)^2}\ll e^{N^\sigma}. $$

Write $\Gamma=H-\text{ diag } h$. For the point-wise estimates, we write 
$$T_N=\tilde T_N+\delta^{2p}\Gamma,$$
$\tilde T_N$ is diagonal. The resolvent equation gives 
$$T_N^{-1}={\tilde T_N}^{-1}+ \delta^{2p} {\tilde T_N}\Gamma {\tilde T_N}^{-1}+ \delta^{4p} {\tilde T_N}\Gamma {\tilde T_N}^{-1}\Gamma {\tilde T_N}^{-1}+ ... $$
If $k\neq k'$, the first term is $0$. Using ($\dag$), $\Gamma$ satisfies 
$$|\Gamma(x, y)|\leq e^{-|x-y|^\gamma}, \, 0<\gamma<1,$$
and $\Vert\Gamma\Vert=\Cal O(1)$. So the above series is norm convergent and moreover 
(2.23) is satisfied. 

Similarly for any fixed 
$$\theta\in (-\infty, -2e^{2(\log N)^2})\cup(2e^{2(\log N)^2}, \infty),$$
since $|\phi|\leq  e^{(\log N)^2}$,
$$|T_\pm|\geq \frac{1}{2}e^{2(\log N)^2},$$
leading to (2.22) and (2.23). 

Note that there is {\it no} excision in $\theta$ or $\phi$, and the proof is direct perturbation.
\smallskip

\noindent Case (c). This is like case (a), and there is no excision in $\theta$ or $\phi$ either.
\smallskip
We observe that in both cases (a) and (c), we expand about the inverse of $\tilde T_N$ ($\{T_\pm(n, j)\}$), 
which may be seen as the Green's function restricted to vertices $(n, j)$, boxes of size $1$.
In all of these boxes, (*) is satisfied. Let us precipitate and say that to treat case (b), we will use 
larger boxes depending on $N$ in the resolvent expansion and that essential will be the upper bounds
on the Green's functions relative to the scale $N$.  
\smallskip 

\noindent Case (b). 
Without loss, one may assume $\theta>0$. 
Since $\theta>0$, $D_+$ as defined in (2.16) satisfies 
$$D_+>\frac{|\theta|+\phi^2}{2}>e^{2(\log N)^2}.$$  
By Schur complement reduction, invertibility of $T_N$ is equivalent to analyze near $0$, 
the effective (symmetric) matrix
$$\aligned \Cal T(\theta, \phi)&=R_-T_NR_- +\{R_-T_NR_+[R_+T_NR_+]^{-1}R_+T_NR_-\}(\theta, \phi)\\
&=\text{diag }[-\theta-n\cdot\omega+M+(j\cdot\lambda+m+\phi)^2]+\delta^{2p}R_-H_NR_-+ \Cal O(\delta^{4p}e^{-2(\log N)^2}),\endaligned\tag 2.26$$
where $H_N$ is the restriction of the convolution matrix $H$ defined in (2.10) (and (2.14)), $R_{\pm}$ are the projections onto the $\pm$ sectors, and the 
last matrix is Lipschitz in $\theta$ and $\phi$ with the stated norm. Below we use a Lipschitz eigenvalue variation technique 
to make norm estimates.

\smallskip
$\bullet$ {\it Norm estimates}

For a fixed $\phi$, one may view (2.26) defining a (weakly) nonlinear eigenvalue problem, and write $\Cal T$ as 
$$\Cal T=-\theta +\Phi(\theta; \phi),$$
where $\Phi$ is self-adjoint and $\Vert \Phi\Vert _{\text{Lip}(\theta)}= \Cal O(\delta^{4p}e^{-2(\log N)^2})$. Diagonalizing $\Phi$ yields Lipschitz 
eigenvalues $E_i=E_i(\theta; \lambda, \omega, m, M; \phi)$, satisfying 
$$\Vert E_i\Vert _{\text{Lip}(\theta)}=\Cal O(\delta^{4p}e^{-2(\log N)^2}),\tag 2.27$$
$$\Vert E_i\Vert _{\text{Lip}(\omega)}=\Cal O(N),\,  \Vert E_i\Vert _{\text{Lip}(M)}=\Cal O(1),\tag 2.28$$
$$\Vert E_i\Vert _{\text{Lip}(\lambda)},  \Vert E_i\Vert _{\text{Lip}(m)}=\Cal O(N|\phi|),\tag 2.29$$
$i\in[-N, N]^4$. (Here we used that the matrix $\Phi$ is Lipschitz in all the parameters.) Using (2.27) and the Lipschitz implicit function theorem \cite{P256, Cl}, this gives a set of Lipschitz functions $\theta_i(\lambda, \omega, m, M; \phi)$, $i\in[-N, N]^4$, 
satisfying the estimates in (2.24), as the $E_i$'s,  such that
if $$|\theta-\theta_i|>2e^{-N^\sigma},$$
for all $i\in[-N, N]^4$, then 
$$\Vert (T_N)^{-1}(\theta)\Vert <e^{N^\sigma}.$$
This gives the measure estimate in (2.20) and the norm estimate in (2.22).

Similarly, for a fixed $\theta$, define $t:=\phi^2$. One may restrict $t$ to $t<10 |\theta|$, as otherwise $\Cal T$ is invertible.   
Write $\Cal T$ as
$$\Cal T=t+\Psi(t),$$
where $\Psi$ is self-adjoint and $\Vert \Psi\Vert _{\text{Lip}(t)}=\Cal O(e^{-(\log N)^2/2})$. Diagonalizing $\Psi$ yields Lipschitz 
eigenvalues $\Cal E_i(t; \lambda, \omega, m, M, \theta)$, satisfying 
$$\Vert \Cal E_i\Vert _{\text{Lip}(t)}=\Cal O(Ne^{-(\log N)^2/2}),\tag 2.30$$
$$\Vert \Cal E_i\Vert _{\text{Lip}(\omega)}=\Cal O(N),\,  \Vert \Cal E_i\Vert _{\text{Lip}(M)}=\Cal O(1),\tag 2.31$$
$$\Vert \Cal E_i\Vert _{\text{Lip}(\lambda)}=\Cal O(N\sqrt {|\theta|}),  \Vert \Cal E_i\Vert _{\text{Lip}(m)}=\Cal O(\sqrt {|\theta|}), \tag 2.32$$
$i\in[-N, N]^4$. As before, using (2.30) and the Lipschitz implicit function theorem yields a set of Lipschitz functions $t_i(\lambda, \omega, m, M; \theta)$, $i\in[-N, N]^4$, 
satisfying the estimates in (2.31) and (2.32) such that
if $$|t-t_i|>2e^{-N^\sigma},$$
for all $i\in[-N, N]^4$, then 
$$\Vert (T_N)^{-1}(\phi)\Vert <e^{N^\sigma},$$
using $\phi=\sqrt t$. This yields the norm estimates in (2.22). Using $d\phi/dt=1/2\phi$ 
gives the measure estimate in (2.21).  One may, evidently, assume $t_i>e^{2(\log N)^2}>0$. Setting $\phi_i=\pm\sqrt{t_i}$, then
$\phi_i$ are Lipschitz satisfying (2.25).

Clearly the above arguments hold for 
$\theta<0$, by interchanging $D_+$ and $D_-$. So we have proven the norm estimates. 
\smallskip

$\bullet$ {\it Pointwise estimates}

To obtain the point-wise estimates in (2.23), we return to the linearized operator $T_N$. 
We use multi-scale induction and cover $[-N, N]^4$ by smaller boxes of size $N_1=N^{1/q}$, for some $q>1$ to be determined.
Assume Lemma~B holds at scale $N_1$. For any fixed $(\theta, \phi)$, 
consider the diagonals of  $T_N(\theta, \phi)$: 
$$T_-(n, j)= -\theta-n\cdot\omega+M+(j\cdot\lambda+m+\phi)^2+\delta^{2p} h,$$
where $h$ is the diagonal of $H_N$ as before. ($|T_+(n, j)|\gg 1$ for $\theta>0$). 
The difference of a pair of diagonals
is
$$T_-(n, j)-T_-(n', j')=(n'-n)\cdot\omega+(j-j')\cdot\lambda[(j+j')\cdot\lambda+2m+2\phi], \, (n, j)\neq (n', j').$$
If $j\neq j'$, then 
$$\aligned |T_-(n, j)-T_-(n', j')|&\geq |\phi| |(j-j')\cdot\lambda|\geq \frac{ \kappa|\phi|}{N^{\rho}}, \quad 0<\kappa<1,\, \rho>3\\
&>e^{(\log N)^2/2}, \endaligned$$
since $\lambda$ is Diophantine and $N\geq |\log\delta|^s\gg 1$, $s>1$. 
So for any fixed $(\theta, \phi)$, there exists a (unique) $\bar j$, such that if
$|T_-(n, j)|<1$, then $j=\bar j$. 

Now for a fixed $\phi$, each $N_1$ box produces $(2N_1+1)^4$ Lipschitz $\theta_i$. 
There are $(2N_1+1)^2$  $j$'s such that 
$$|j-\bar j|\leq N_1.\tag 2.33$$
This leads to $(2N_1+1)^6$ Lipschitz $\theta_i$ for cubes $\Lambda_{(0, j)}$ centered 
at $(0, j)$ of size $N_1$, with $j$ satisfying (2.33). 

Let $\Lambda_{(0, j)}$ be such a cube, then its vertical translates
$$\Lambda_{(n, j)}=(n, 0)+\Lambda_{(0, j)},$$
have corresponding Lipschitz functions 
$$\theta_i^{(n)}=n\cdot\omega+\theta_i.$$
For $n, n'\in[-N, N]^2$, $n\neq n'$,
$$|\theta_i^{(n)}-\theta_i^{(n')}|= |(n-n')\cdot\omega|\geq \frac {\kappa}{N^{\rho}}\gg e^{-N_1^{\sigma'}},\, 0<\kappa<1,\, \rho>3,\, \sigma'>0,$$
since $N\gg 1$ and $N_1\gg 1$. ($\sigma'$ will be taken to be $\sigma'=\sigma/q$, see below.)

So for any fixed $(\theta, \phi)$, there can be at most $(2N_1+1)^6$ pair-wise
disjoint $N_1$ boxes $\Lambda$ such that 
there exists a Lipschitz function $\tilde\theta$ for $\Lambda$,
$$\tilde\theta\in\{\theta_i^{(n)}, |n|\leq N, |i|\leq (2N_1+1)^6\}:=S,$$
and 
$$|\theta-\tilde\theta|<e^{-N_1^{\sigma'}}.$$
Taking $q\geq 7$, this yields a family $\Cal F$ of at most 
$N^{6/7}$ pair-wise disjoint bad $N_1$ boxes $\Lambda$ in $[-N, N]^4$. 
If $\Lambda\notin \Cal F$, (2.22) and (2.23) are satisfied at scale $N_1$.

Now fix $\phi$,  and $\theta$ such that 
the norm estimates (2.22) holds for all scales in $[N_1, N]$
in coverings for $[-N, N]^4$, by taking $\sigma'=\sigma/q$, 
Theorem~2.1 \cite{L} applies, leading to (2.23). 

Similarly for fixed $\theta$, 
for  $\phi$ such that the norm estimates (2.22) holds for all scales in $[N_1, N]$,
Theorem~2.1 \cite{L} gives (2.23). Iteration then completes the proof. (Cf. \cite{L}.)

(Note that here we used also that if $\theta$ and $\phi$ satisfy (b) for some $N=N_0$, then they satisfy (b) for all $N>N_0$.)



\hfill $\square$
\enddemo
\smallskip
\noindent{\it Remark.} 
Note that Lemmas~A and B are both stated with the operator linearized at $u_N$ (and $v_N$). Using
($\dag\dag$), one may, however, substitute $u_N$ by any 
$u_{N'}$ for all $N'\geq N$. That ($\dag\dag$) holds will 
be proven in sect.~3 using the double exponential convergence of the Newton scheme. 
\bigskip
\subheading {2.5 The $(\lambda, \omega)$-coordinates}

For invertibility in region (i), we need to change strategy and use 
semi-algebraic geometry. For this purpose, the relation of $(\lambda, \omega)$ 
and $(\lambda, m, M)$ should be algebraic, and we {\it may no longer} work in
$(\lambda, \omega, m, M)\in\Bbb R^6$, as in sects.~2.3 and 2.4. 
For a given $u$ and $v$, we view the $Q$-equations:
$$\omega_k=(h_k\cdot\lambda+m)^2+M+\delta^{2p}\frac{[{(u*v})^{*p}*u] (\lambda, m, M)}{a_k}(-e_k, h_k),\, k=1, 2,\tag 2.34$$
(together with the identity map in $\lambda$) as defining a map from $(\lambda, m, M)$ to $(\lambda, \omega)$.
It will be convenient to work in the $(\lambda, \omega)$-coordinates.
The final estimates, however, will be in terms of $(\lambda, m, M)$.  
\bigskip
Let us start with a definition.
\smallskip
\noindent {\bf Definition.} We call a map $\mu$ from an open set $I\in\Bbb R^d$
to an open set $J\in\Bbb R^d$, {\it bi-rational}, if both $\mu$ and $\mu^{-1}$ are rational. 

For example, it follows from the $Q$-equation in (2.34)
that $\omega$ is rational in $(\lambda, m, M)$, when $u=u^{(0)}$ and $v=u^{(0)}$. Furthermore,
it may be readily solved for $m$ and $M$, yielding
$$\aligned m(\lambda, \omega) &= \frac{\Omega_1-\Omega_2-(h_1\cdot\lambda)^2+(h_2\cdot\lambda)^2}{2(h_1-h_2)\cdot\lambda},\\
M(\lambda, \omega)&=\frac{\Omega_1+\Omega_2}{2}-\frac{[(h_1-h_2)\cdot\lambda]^2}{4}-\frac{(\Omega_1-\Omega_2)^2}{4[(h_1-h_2)\cdot\lambda]^2},\endaligned\tag 2.35$$
where $$\Omega_k=\omega_k-\delta^{2p}\frac{{(u^{(0)}*v^{(0)})^{*p}*u^{(0)}}}{a_k}(-e_k, h_k),\, k=1, 2,\tag 2.36$$
provided $$(h_1-h_2)\cdot\lambda\neq 0.$$
So when $u=u^{(0)}$ and $v=v^{(0)}$, 
$$ (\lambda, m, M)\leftrightarrow (\lambda, \omega),$$
is a bi-rational map. 

Note that the term {\it bi-rational} isomorphism generally pertains to two Zarisky open sets, cf. p26 \cite{Ha}. Here 
we borrow the term and use it for open sets on $\Bbb R^d$. 
To continue, we do not use the relation (2.34), instead we make assumptions, 
{\it separately} on properties of changes of variables:
$$(\lambda, m, M){\leftrightarrow}  (\lambda, \omega),$$
and on $u_N$ and $v_N$. 
Subsequently, in the nonlinear analysis in sect.~3, we verify these assumptions using (2.34).

For simplicity of notation, we use the same letter to
denote a function (operator) in $(\lambda, m, M)$ or $(\lambda, \omega)$ variables.
\smallskip
\noindent{\it Remark.} We note that the diffeomorphism:  
$$ (\lambda, m, M)\leftrightarrow (\lambda, \omega),$$
involves {\it all} $4$ variables, and is {\it not} of the form $\omega=\omega(m, M)$,
for example. This is contrary to previous cases in e.g., \cite{BGS}\cite{B2}\cite{BK},
and is the reason why we make bi-rational approximations.
\bigskip
\subheading {2.6  The Main Lemma and region (i)}
\smallskip
Write $k=(n, j)\in\Bbb Z^2\times\Bbb Z^2:=\Bbb Z^4$. Below we call connected open sets on $\Bbb R^d$, {\it intervals}. 
Let $I$ be an interval in $\Bbb R^4$. To proceed, we define two properties:
\bigskip 
\noindent {\bf {Property B:}}

We say that a sequence of {\it bi-rational} maps $\Cal M_N$: 
$$\aligned (\lambda, m, M) &\overset {\Cal M_N}\to { \underset {\Cal M_N^{-1}} \to\longleftrightarrow } (\lambda, \omega),\\
I&\leftrightarrow\Cal I,\endaligned\tag 2.37$$
has Property B, if 
$$\align&\det \big [\frac{\partial (\lambda, \omega)}{\partial(\lambda, m, M)}\big]=\Cal O(1),\tag 2.38\\
&\det \big [\frac{\partial (\lambda, m, M)}{\partial(\lambda, \omega)}\big]=\Cal O(1),\tag 2.39\endalign$$
uniformly in $N$;
$$|\Cal M_N-\Cal M_{N+1}|_\infty\leq e^{-\tilde\gamma N},\tag 2.40$$
where $\tilde\gamma>0$, as in Property V1; and 
$$\aligned &\text{deg }m(\lambda, \omega)\lesssim e^{(\log N)^3},\\
&\text{deg }M(\lambda, \omega)\lesssim e^{(\log N)^3}.\endaligned\tag 2.41$$
\bigskip
\noindent {\bf {Property V3:}}

We say that a sequence of real {\it rational} functions: 
$$w_N:\, I\mapsto \ell^2(\Bbb Z^{4}),\, N\in\Bbb N,$$
has Property V3, if it satisfies Property V1 on $I$, and 
$$\text{deg }w_N\circ \Cal M^{-1}_N\lesssim e^{(\log N)^3},\tag 2.42$$
provided the map $\Cal M_N$ has Property B. 
\bigskip
\proclaim{Main Lemma} Let $0<\delta\ll 1$, and $I$ be an interval in $(\lambda, m, M)$. Assume that on $I$, there is a sequence of bi-rational maps
$$\align (\lambda, m, M) &\overset {\Cal M_N}\to { \underset {\Cal M_N^{-1}} \to\longleftrightarrow } (\lambda, \omega),\\
I&\leftrightarrow\Cal I,\endalign$$
satisfying Property B;
and two sequences of real rational functions $u_N$ and $v_N$:
$$u_N:\, I\mapsto \ell^2(\Bbb Z^{4}),$$
$$v_N:\, I\mapsto \ell^2(\Bbb Z^{4}),$$
satisfying Property V3. 
Then for all $N\geq N_0$, defined as in (2.19),  there exists $\Cal G_N\subset\Cal I$, a semi-algebraic set in $(\lambda, \omega)$, satisfying
$$\aligned \text{deg }&\Cal G_N\leq e^{(\log N)^3}, \\
\text{meas }(\Cal G_{N-1}\backslash& \Cal G_N)\leq N^{-2c},\, c>0.\endaligned\tag 2.43$$
Let $$I\supset \Gamma_{N}=\Cal M_N^{-1}(\Cal G_N),$$
and $$I\supset \Gamma'_{N-1}=\Cal M_{N}^{-1}(\Cal G_{N-1}),$$
$$\text{meas }(\Gamma'_{N-1}\backslash \Gamma_N)\leq N^{-c},c>0.\tag 2.44$$
For any $(\lambda, \omega)\in\Cal G_N$, any $(\lambda, m, M)\in \Gamma_N$, there exists a subset $\Theta_N\subset \Bbb R^2$, whose sectional measures
satisfy
$$\text{meas } [\theta| \text{fixed }\phi, (\theta, \phi)\in \Theta_N] \leq e^{-N^\tau}\, (\tau>0),\tag 2.45$$
and $$\text{meas } [\phi| \text{fixed }\theta, (\theta, \phi)\in \Theta_N]\leq e^{-N^\tau} \, (\tau>0).\tag 2.46$$
If $(\theta,\phi)\notin \Theta_N$, the linearized operator $F':=T$ in (2.8)-(2.10), after truncation and with
$u=u_N$ and $v=v_N$, satisfy the estimates
$$\Vert [T_N(\theta,\phi)]^{-1}\Vert_{\text{Op}}\leq e^{N^\sigma}\, (0<\tau<\sigma<1),\tag 2.47$$
$$|[T_N(\theta,\phi)]^{-1}(k,k')|\leq e^{-|k-k'|^\gamma}\, (0<\tau<\sigma<\gamma<1), \forall |k-k'|>N/10,\tag 2.48$$
where $T_N$ stands for either $T_N(\lambda,\omega; u_N, v_N)$ or $T_N(\lambda, m, M; u_N, v_N)$.
\endproclaim
\bigskip
When applying the Main Lemma in sect.~3.5, we shall fulfil the assumptions on 
$\Cal M_N$ $u_N$ and $v_N$, namely Properties B and V3, by using the structure 
presented by (2.11)-(2.13) and (2.35)-(2.36), together with the double exponential convergence of the Newton scheme.

More precisely, it will follow from (2.11) and (2.12) that for all $N$, on {\it good} intervals $I$, to be constructed in sect.~3, the map $\Cal M_N$ defined in (2.13), equivalently (2.34):
$$(\lambda, m, M)\overset { \Cal M_N(u_N,  v_N) }\to {\longmapsto}  (\lambda, \omega),$$ 
is {\it rational}. Moreover, using another simple Newton scheme (with no small-divisors) to solve (2.13), equivalently (2.34),
we will show in sect.~3.4 that 
on the same good intervals, 
the diffeomorphism $\Cal M_N$ admits an inverse $\Cal M_N^{-1}$,
which can be {\it well approximated} by a rational map, 
and that it suffices for our application, cf. Lemma~3.3.

\noindent{\it Remark.} 
It is fairly standard that (2.34) defines a diffeomorphism on $(\lambda, m, M)$, by using implicit (inverse) function theorem, see sect.~3.

In view of Lemmas A and B, we only need to deal with region (i).
The difference of the Main Lemma in region (i) with that in \cite{W3} is the 
assumption 
of a family of bi-rational changes of variables: 
$$(\lambda, m, M) \overset {\Cal M_N}\to { \underset {\Cal M_N^{-1}} \to\longleftrightarrow } (\lambda, \omega).$$ 
The final estimates for each $N$ are, however, always expressed in the
original parameters $(\lambda, m, M)$. 

In region (i), the proof relies on tools developed in the course of Proposition~2.2 in \cite {B3}.
Even though we do not repeat the proof, 
we give, nonetheless, the definition of semi-algebraic sets and state the basic algebraic lemmas 
used there.  
\bigskip
\subheading {2.7  Semi-algebraic sets}
\smallskip
\noindent{\bf Definition.} A set $S$ is called semi-algebraic if it is a finite union of sets defined by a 
finite number of polynomial equalities and inequalities. More specifically, let $\Cal P=\{P_1, P_2, ..., P_s\}\subset \Bbb R[x_1, x_2, ..., x_n]$
be a family of $s$ real polynomials of degree bounded by $\kappa$. A (closed) semi-algebraic set $S$ is given by an expression
$$S=\bigcup_j \bigcap_{\ell\in\Cal L_j} \{P_\ell s_{jl}0\},\tag S$$
where $\Cal L_j\subset \{1, 2, ..., s\}$ and $s_{jl}\in\{\geq, =,\leq\}$ are arbitrary. We say that $S$ as introduced above has degree at most $s\kappa$ and
its degree $B$ is the minimum $s\kappa$ over all representations (S) of $S$. 

The following is a special case of Theorem 1 in \cite{Ba}, cf. Theorem 9.3 in Chap.~9 \cite{B2}.

\proclaim{Lemma~2.1}Let $S\subset \Bbb R^n$ be as in {\rm (S)}. Then the number of connected components of $S$ does not exceed 
$\Cal O(s\kappa)^n$. 
\endproclaim

The two properties of semi-algebraic sets that play a central role here are the Tarski-Seidenberg principle, which states that the projection 
of a semi-algebraic set of $\Bbb R^n$ onto $\Bbb R^{n-1}$ is semi-algebraic; and the Yomdin-Gromov triangulation theorem of these sets. 
They are both stated in \cite{B3}, cf. the references therein. (For the complete proof of the Yomdin-Gromov triangulation theorem, see \cite{BiN}, cf., also the earlier 
paper \cite{Bu}.) We do not repeat them here, except the following consequences for {\it thin} sets.

\proclaim{Lemma~2.2}
Let $S\subset [0, 1]^{n_1}\times [0, 1]^{n_2}:=[0, 1]^{n}$, be a semi-algebraic set of degree $B$ and $\text{\rm meas}_{n} S <\eta, \log B\ll
\log 1/\eta$.
Denote by $(x, y)\in [0, 1]^{n_1}\times [0, 1]^{n_2}$ the product variable.
Fix $\epsilon>\eta^{1/n}$.
Then there is a decomposition
$$
S =S_1 \bigcup S_2,
$$
with $S_1$ satisfying
$$
\text{\rm meas}_{n_1}(\text{Proj}_x S_1)<B^K\epsilon \quad (K>0),
$$
and $S_2$ the transversality property
$$
\text{\rm meas}_{n_2}(S_2\cap L)< B^K\epsilon^{-1} \eta^{1/n}\quad (K>0),
$$
for any $n_2$-dimensional hyperplane $L$ in  $[0, 1]^{n_1+n_2}$ such that
$$
\max_{1\leq j\leq n_1}|\text{Proj}_L (e_j)|< \frac 1{100}\epsilon,
$$
where $e_j$ are the basis vectors for the $x$-coordinates.
\endproclaim

Lemma~2.2 is the basic tool,
underlining the semi-algebraic techniques used in the subject. 
It is stated as (1.5) in \cite{B3}, cf., Lemma~9.9 \cite {B2} and
Proposition ~5.1 \cite{BGS} and their proofs. The $\epsilon$  
measures {\it steepness}. If $x$ is the horizontal direction
and $y$ vertical, then the set $S_1$ can be viewed as the 
{\it vertical component}, and $S_2$, {\it horizontal}. In applications
$\eta\ll \epsilon$, so $S_1$ is usually the dominant term, leading 
to much bigger loss of measure. 

In our application, $x$ will be 
identified with $(\lambda, \omega)$, and $y$, $(\theta, \phi)$; 
$S$ with the complement of $\Cal G_N$.
We shall use the covariance (2.17), so $\epsilon\sim |(n, j)|^{-1}$,
leading to the steepness notion mentioned in the Introduction. 
However, this covariance needs to be applied to {\it any}
starting point $(\theta_0, \phi_0)$. Clearly $S$ 
depends on $(\theta_0, \phi_0)$. So we may not
apply Lemma~2.2 as it is, which would lead to 
sets $S_1$ and $S_2$ {\it dependent} on $(\theta_0, \phi_0)$. We need
to {\it eliminate} the variable $(\theta_0, \phi_0)$. This is achieved 
by the two lemmas below. 
\bigskip
\proclaim{Lemma~2.3}
Let $A\subset [0, 1]^{n+r}$ be semi-algebraic of degree $B$ and such that  
$$ \text{for each }t\in[0, 1]^r, \, \text{meas}_nA(\cdot\,,t)<\eta, \eta>0.$$
Then $$\Cal A:=\{(x_1, x_2, ..., x_{2^r})|A(x_1)\cap...\cap A(x_{2^r})\neq \emptyset\}\subset[0, 1]^{n2^r},$$
is semi-algebraic of degree at most $B^C$ and measure at most 
$$\eta_r=B^C\eta^{n^{-r}2^{-\frac{r(r-1)}{2}}},$$
with $C=C(r)>1$. 
\endproclaim

\noindent Lemma~2.3, stated as Lemma~1.18, and proven, in \cite{B3}, is a variable reduction lemma, eliminating the $r$-dimensional variable $t$.
It is worth noting that $2^r$ copies of $A$ are used. The measure is, however, in
$n2^r$ dimensions; while we need the measure of a $n$-dimensional section of $\Cal A$.
Lemma~1.20 in \cite{B3} (and proven there) serves this purpose, and is stated below.

\proclaim{Lemma 2.4} Let $A\subset [0, 1]^{nd}$ be a semi-algebraic set of degree $B$ and 
$$\text{meas}_{nd}A<\eta.$$
Let $w_i\in [0,1]$, $i=1, 2, ..., n$,
and $$w=(w_1, w_2, ..., w_n)\in [0, 1]^{n}.$$
Let $k_i\in \Bbb Z$, $i=1, 2, ..., n$,
and $$k=(k_1, k_2, ..., k_n)\in \Bbb Z^{n}.$$
Denote by $\{\cdot\}$, the fractional part of a real number in $[0, 1)$, and   
$$kw:=(\{k_1w_1\}, \{k_2w_2\}, ..., \{k_nw_n\}).$$
Let $K_1$, $K_2$, ..., $K_{d-1}\subset \Bbb Z^{n}$ be finite sets with the following properties:
$$\min_{1\leq \ell \leq n}|k_\ell|>[B \max_{1\leq \ell' \leq n}|m_{\ell'}|]^C,$$
if $k\in K_i$ and $m\in K_{i-1}$, $i=2, ..., d-1$, and where $C=C(n, d)$. Assume also 
$$\frac{1}{\eta}>\max_{k\in K_{d-1}}|k|^C.$$
Then 
$$\aligned &\text{meas }\{w\in [0, 1]^{n}|(w, k^{(1)}w, ..., k^{(d-1)}w)\in A \text{ for some } k^{(i)}\in K_i\}\\
<&B^C\delta,\endaligned$$
where $$\frac{1}{\delta}=\min_{k\in K_1}\min_{1\leq \ell\leq n}|k_\ell|.$$
\endproclaim
\smallskip
\noindent{\it Remark.} Semi-algebraic variable elimination is the main development in \cite{B3}, compared 
to the previous semi-algebraic techniques in, e.g., \cite{BGS} and 
Chap.~18-20 \cite{B2}.   
\bigskip
\subheading {2.8  Proof of the Main Lemma}
\smallskip
\demo{Proof of the Main Lemma}
For the initial scales defined in (2.18), we may invoke Lemma A. 
\smallskip
For larger scales $N$, $N\geq |\log\delta|^s$ $(s>1)$, there are two cases, as mentioned earlier: 
\item {(i)} $|\theta|\leq 2e^{2(\log N)^2}$ and $|\phi|\leq 2e^{(\log N)^2}$;
\item{(ii)} otherwise.

After restricting $\omega$ to be the $\omega$ component of the map $\Cal M_N$, i.e., 
$$\omega=\Cal M_{N,\omega}(\lambda, m, M),$$
we use Lemma~B to address case (ii). The Diophantine conditions lead to
$(2N+1)^4-2$ monomials in $\lambda$ or $\omega$, and excision
of a set $\Cal S_N$ in $\Cal I$ satisfying
$$\text{meas } (\Cal S_N\backslash\Cal S_{N-1})\leq \frac{\kappa}{N}.\tag 2.49$$
Let 
$$\aligned S_N&=\Cal M_N^{-1}(\Cal S_N),\\
S'_{N-1}&=\Cal M_N^{-1}(\Cal S_{N-1}),\endaligned$$ 
then (2.39) yields
$$\text{meas } (S_N\backslash S'_{N-1})\leq \frac{\Cal O(1)}{N}.$$

On $\Cal I\backslash \Cal S_N$ and $I\backslash S_N$, application of Lemma~B yields
(2.45)-(2.48). We are now only left with case (i). 

As in \cite{W3},
we do not use the special structure of the $\Bbb Z^4$ action in (2.17) and double the dimension of $(\theta, \phi)\in\Bbb R^2$ 
to $(\tilde\theta, \tilde\phi)\in\Bbb R^4$:
$$\aligned n\cdot\omega&=n_1\omega_1+n_2\omega_2\mapsto n_1\omega_1+\tilde\theta_1+n_2\omega_2+\tilde\theta_2;\\
j\cdot\lambda&=j_1\lambda_1+j_2\lambda_2\mapsto j_1\lambda_1+\tilde\phi_1+j_2\lambda_2+\tilde\phi_2.\endaligned$$
It follows that 
$$\aligned \theta&=\tilde\theta_1+\tilde\theta_2,\\
 \phi&=\tilde\phi_1+\tilde\phi_2.\endaligned$$
The sectional measures in $(\theta, \phi)$ can be derived from that in 
 $$(\tilde\theta, \tilde\phi)=(\tilde\theta_1, \tilde\theta_2, \tilde\phi_1, \tilde\phi_2).$$
 This can be seen as below.

 Let $S\subset \Bbb R^4$, assume, for example,
 $$\text{meas }S(\tilde \theta_1|\text{ fixed } \tilde\theta_2, \tilde\phi_1, \tilde\phi_2)\leq \zeta,$$
 then 
$$\aligned \text{meas }S(\tilde \theta_1+\tilde \theta_2|\text{ fixed } \tilde\phi_1, \tilde\phi_2)&=\text{meas }S(\theta|\text{ fixed } \phi)\\
&\leq \zeta,\endaligned$$
using that the set $S$ is independent of $\tilde\theta_1-\tilde \theta_2$.   Clearly similar argument leads to the estimate in $\phi$. 
 
To prepare for the proof in region (i), we therefore set:
 $$\tilde\theta_{s,i}:=\theta_s-\tilde\theta_{i'},\, i,i'=1,2,\, i\neq i',$$
 and
 $$\tilde\phi_{s,i}:=\phi_s-\tilde\phi_{i'},\, i, i'=1,2,\, i\neq i',$$
 where $\theta_s$ and $\phi_s$ are the Lipschitz functions 
 constructed in Lemma~B. After the doubling of dimensions, the analogue of Lemma~B
 in the $(\tilde\theta, \tilde\phi)$ variables clearly holds. 
 \bigskip
 
Returning to region (i), we first note that if on $\Cal G_{N}$, 
$$\aligned \Vert [T_{N}(\lambda, \omega; u_{N}, v_{N})(\theta, \phi)]^{-1}\Vert_{\text{Op}} &\leq e^{N^\sigma},\,\qquad\,  (0<\tau<\sigma<1),\\
|[T_{N}(\lambda, \omega; {u_{N}, v_{N}})(\theta, \phi)]^{-1}(k, k')| &\leq e^{-|k-k'|^\gamma},\, (0<\tau<\sigma<\gamma<1),\,\forall |k-k'|>N/10,\endaligned\tag 2.50$$
where $$(\lambda, \omega)=\Cal M_{N}(\lambda, m, M),$$ is the image under the $N$th bi-rational map, then using Properties V3 and B (2.40),   
we have 
$$\aligned \Vert [T_{N}(\lambda, \omega; u_{N'}, v_{N'})(\theta, \phi)]^{-1}\Vert_{\text{Op}} &\leq (1+e^{-\bar\gamma N}) e^{{N}^\sigma},\\
|[T_{N}(\lambda, \omega; {u_{N'}, v_{N'}})(\theta, \phi)]^{-1}(k,k')| &\leq (1+e^{-\bar\gamma N})e^{-|k-k'|^\gamma},\,\forall |k-k'|>N/10,\endaligned\tag 2.51$$
for some $\bar\gamma$, $0<\bar\gamma<\tilde\gamma$,
with $$(\lambda, \omega)=\Cal M_{N'}(\lambda, m, M),$$ for all $N'>N$. So 
one may view $\Cal G_N$ as a {\it fixed} set (with respect to the maps $\Cal M_{N'}$, $N'\geq N$),
after replacing (2.50) by (2.51).

In the $(\lambda, \omega)$-coordinates, we are now in the same setting as the Main Lemma in \cite{W3}.  The Lemma and the proof 
there apply. Recall that the proof uses three scales: 
$N_1$, $N_2=N_1^{2/\tau}$ and $N_3=e^{N_1^\tau}$. Precisely, assume (2.47) and (2.48) hold at scales $N_1$ on the set $\Cal G_{N_1}$, and $N_2$, on 
$\Cal G_{N_2}$. In order that (2.47) and (2.48) hold at all scales $N\in[N_3, N_3^2]$, the excision is made on the set $\Cal G_{N_1}$. Call the resulting 
set $\Omega_3$:
$$\text{meas }\Cal G_{N_1}\backslash\Omega_3\leq N_3^{-c'},0<c'<1.\tag 2.52$$

Set $\Cal G'_N=\Omega_3\cap \Cal G_{N_2}$, so that (2.47) and (2.48) are available at scale $N_2$.   
Remark that $\Cal G'_N$ is obtained from $\Cal G_{N_1}$,  and 
$$N_1\sim (\log N)^{1/\tau}\ll N.$$ 
Set $$\Cal G''_N=\Cal G'_N\cap (\Cal I\backslash \Cal S_N),$$
so that (2.45)-(2.48) hold on region (ii).  For later nonlinear analysis, it will be
convenient to have nested sets. So we define
$$\Cal G_N=\Cal G''_N\cap_{N'<N}\Cal G_{N'}.$$

Since there are $\Cal O(N^4)\ll e^{(\log N)^3}$ Diophantine monomial conditions, and from (2.49)
$N^{-1}\ll N^{-c'}$, 
(2.43) is satisfied, with $c=16\tau$.  On $\Cal G_N$, (2.45)-(2.48) hold at scale $N$ in the entire $(\theta, \phi)$ space, and there is the nested property:
$$\Cal I\supseteq \cdots\supseteq\Cal G_{N_1}\cdots\supseteq\Cal G_{N_2}\cdots\supseteq \Cal G_{N-1}\supseteq\Cal G_{N}\supseteq \cdots 
\supseteq \Cal G_{N+\ell}\supseteq \cdots, \quad \ell\in\Bbb N.\tag 2.53$$
The bound on the degree follows from the definition above (2.53)
and the observation that each intersection leads to at most a polynomial factor in $N'$.

Using (2.39) and (2.43), yields (2.44), provided $1>s\tau>1/16$. 
From (2.38)-(2.40):
$$|\Cal M^{-1}_{N-1}-\Cal M^{-1}_N|_\infty\leq e^{-\tilde\gamma (N-1)},\tag 2.54$$
since $N\geq |\log\delta|^s$, $s>1$.
Using (2.54) and the same argument leading (2.50) to (2.51), but in the $(\lambda, m, M)$ variables, yield that if  
$$\aligned \Vert [T_{N}(\lambda, m, M; u_{N}, v_{N})(\theta, \phi)]^{-1}\Vert_{\text{Op}} &\leq e^{N^\sigma},\,\qquad\,  (0<\tau<\sigma<1),\\
|[T_{N}(\lambda, m, M; {u_{N}, v_{N}})(\theta, \phi)]^{-1}(j,j')| &\leq e^{-|j-j'|^\gamma},\, (0<\tau<\sigma<\gamma<1),\,\forall |j-j'|>N/10,\endaligned$$
then
$$\aligned \Vert [T_{N}(\lambda, m, M; u_{N'}, v_{N'})(\theta, \phi)]^{-1}\Vert_{\text{Op}} &\leq (1+e^{-\bar \gamma N}) e^{{N}^\sigma},\\
|[T_{N}(\lambda, m, M; {u_{N'}, v_{N'}})(\theta, \phi)]^{-1}(j,j')| &\leq (1+e^{-\bar\gamma N})e^{-|j-j'|^\gamma},\,\forall |j-j'|>N/10,\endaligned$$
for some $\bar\gamma>0$, and all $N'>N$.

This shows that one may redefine the family of $\Gamma_N$ 
so that 
$$I\supseteq \cdots\supseteq \Gamma_{N-1}\supseteq \Gamma_{N}\supseteq\cdots \supseteq \Gamma_{N+\ell}\supseteq \cdots, \quad \ell\in\Bbb N,\tag 2.55$$
satisfying $$\text{meas }(\Gamma_{N-1}\backslash \Gamma_N)\leq 2N^{-c},c=16\tau>0,\tag 2.56$$
and on $\Gamma_N$, (2.47) and (2.48) are essentially preserved. This completes the proof, after lowering $\sigma$ and $\beta$, if needed.
\bigskip
Finally note that in terms of Lemma~2.4, we are here in the setting $n=4$ and $d=2^4+1=17$. 
\hfill $\square$
\enddemo
\bigskip
Let us encapsulate the observations in (2.53) and (2.55), which foreshadow the nonlinear analysis in sect.~3. 

\proclaim{Corollary to the Main Lemma}
The sets $\Cal G_N$ and $\Gamma_N$ in the Main Lemma, may be chosen to exhibit the nested properties in
(2.53) and (2.55) respectively. 
\endproclaim
\bigskip
\subheading{2.9 Comparison of the proofs in regions (i) and (ii, b)}

It is instructive to find a common thread for the seemingly disparate proofs in these two cases. In (i), broadly speaking, Lemma~2.4 is used to
control the number of bad regions of smaller scales $N'\ll N$, by constructing the set $A$ from $\Theta_N$, setting $w=(\lambda, \omega)$ and excise in $(\lambda, \omega)$. 
It should be noted that in (i), we have very little information on $\Theta_N$, and consequently $A$, aside 
from it being semi-algebraic and of small measure. This is in contrast with (ii, b), where we have good control of the bad set
$\Theta_N$, through the families of Lipschitz functions $\theta_i$ and $\phi_i$. Call $(\lambda, \omega)$ the horizontal variables, and $(\tilde \theta, \tilde \phi)$ vertical. 
Putting into the context of Lemma~2.2, one may view
$\theta_i$ and $\phi_i$ as manifesting, in some sense, the {\it absence} of the vertical set $S_1$, which usually dominates 
the (loss of) measure. 

\noindent{\it Remark 1.} 
In case (i), bi-rationality is only used in the algebraic part 
of the proof, to deduce that for {\it any fixed} $(\tilde \theta, \tilde \phi)$, there are only sublinear bad appropriately smaller boxes. This, however, permits a perturbation
of size $\Cal O(e^{-N_1})$ here, in the estimates for the good boxes, cf. the paragraph containing (2.52). (This is also used in the proof of the Main Lemma in \cite{W3}.)

\noindent{\it Remark 2.} It is useful to note that in the multi-scale proofs, the expansion of scales is always polynomial or sub-exponential.
This (at most) sub-linear rate is  {\it canonical} to multi-scale analysis. It applies 
as well to the nonlinear analysis, cf. the scales in sect.~3.5, in particular Lemma~3.4. 

\head{\bf 3. Nonlinear analysis -- proof of the Theorem }\endhead

We can now construct space-time quasi-periodic solutions of the form (A) to the
NLS in (1.1). From sect.~2.1, this leads to solve the nonlinear 
matrix equation in (2.4), which we denote by $F(u, v)=0$ in (2.6)-(2.7). 
(Recall that $u$ and $v$ now stand for $\hat u$ and $\hat v$.) 
The $P$-equations in (2.6) are used to solve for the amplitudes; while the $Q$-equations
in (2.7), the frequencies, as in (2.13). 

We use the Newton scheme
in (2.11)-(2.12) to solve the $P$-equations. Let $A$ be a large positive integer, and 
$A^r$, $r=1, 2, ...$, the geometric sequence of scales.  
Denote by $u^{(r)}$ the $r$th approximation, and the increment 
$$\aligned \Delta u^{(r)}&=u^{(r)}-u^{(r-1)},\, r\geq 1,\\
 \Delta v^{(r)}&=v^{(r)}-u^{(r-1)}, \, r\geq 1.\endaligned\tag N1$$
 Use (2.13) and set  
$$\omega^{(r)}_k=(h_k\cdot\lambda+m)^2+M+\delta^{2p}\frac{[{(u^{(r-1)}*v^{(r-1)}})^{*p}*u^{(r-1)}] (\lambda, m, M)}{a_k}(-e_k, h_k),\, k=1, 2, \, r\geq 1.\tag N2$$ 
We define 
$$\Delta \pmatrix u^{(r)}\\v^{(r)}\endpmatrix= -[F'_N(u^{(r-1)}, v^{(r-1)}, \omega^{(r)})]^{-1}F(u^{(r-1)}, v^{(r-1)}, \omega^{(r)}),\tag N3$$
where $N=A^{r}$ and $F'_N(u^{(r-1)}, v^{(r-1)}, \omega^{(r)})$ is the restriction of the linearized operator
$F'(u^{(r-1)}, v^{(r-1)}, \omega^{(r)})$
to the cube $[-N, N]^{4}\times\{0, 1\}$.
\smallskip
Equations (N1)-(N3), together with $u^{(0)}$, $v^{(0)}$ and $\omega^{(0)}$, iteratively solve the $Q$ and the $P$-equations,
provided (N3) is well-defined for all $r$, and the resulting series is absolutely convergent. From (N3), the invertibility of $F'_N$
is central, and indeed will be provided by the Main Lemma, which leads to double exponential convergence.  
\smallskip
By inverse (implicit) function theorem, the solutions to the $Q$-equations in (N2) 
yield diffeomorphisms, as in previous works, for example, Chap.~18 \cite{B2} and \cite{W3}.
This, however, does not suffice for 
applications of the Main Lemma here, which require an algebraic setting.
We henceforth use another Newton scheme to make bi-rational 
approximations to the maps defined by (N2). This is a new development
of the method.
\bigskip
\subheading{3.1 The induction hypothesis}

Recall that $(\lambda, m, M)=(\lambda_1, \lambda_2, m, M)\in (0, 2\pi)^4$. We call $(0, 2\pi)^4$,
the {\it entire space}. Since $u^{(r)}$ and $v^{(r)}$ 
 are complex conjugates, to save space, below we will generally make statements on $u^{(r)}$ only.

Let $A$ be a large positive integer. The nonlinear analysis consists in showing that the following are satisfied for all $r\geq 1$ and fixed small $\delta$: 

On the entire $(\lambda, m, M)$ space:
\item{(Hi)} $\text{supp }u^{(r)}\subseteq B(0, A^r)$ ($\text{supp }u^{(0)}\subset B(0, A)$).
\item{(Hii)} $\Vert \Delta u^{(r)}\Vert<\delta_r$,
$\Vert \partial \Delta u^{(r)}\Vert<\tilde\delta_r$ with $\delta_{r+1}\ll\delta_{r}$ and $\tilde\delta_{r+1}\ll\tilde\delta_{r}$,
where $\partial$ denotes $\partial_x$, $x$ stands for $\lambda_1, \lambda_2, m$ or $M$,
and $\Vert\,\Vert:=\sup_{(\lambda, m, M)}\Vert\,\Vert_{\ell^2(\Bbb Z^{4})}$.
\item{(Hiii)} $|u^{(r)}(k)|<\delta e^{-|k|^\gamma}\, (0<\gamma<1)$. 
\bigskip
Use (Hi-iii) and
view (N2) as defining the $\omega$ component of a map $\Cal M_r$: 
$(\lambda, m, M)\mapsto (\lambda, \omega)$;
the $\lambda$ component being the identity map. 
Set 
$$\omega^{(r)}:=\Cal M_{r, \omega} (\lambda, m, M), \, r\geq 1,\tag 3.1$$
where $ \Cal M_{r, .}$ denote components of the map. ($\omega^{(0)}$ as defined in (2.2).)
Then $\omega^{(r)}$ is $C^1$ in $(\lambda, m, M)$ on the entire space. 
Moreover  by (Hii), 
$$| \omega^{(r)}-\omega ^{(r-1)} |\lesssim \Vert u^{(r-1)}-u^{(r-2)}\Vert <\delta_{r-1},\, r\geq 1,\tag 3.2$$
where $u^{(-1)}$ is set to be: $u^{(-1)}\equiv 0$.
\smallskip
Below we continue with the assumptions on the {\it restricted} intervals in $(\lambda, m, M)$, where one could construct approximate solutions.

\item{(Hiv)} There is a set $\Lambda_r$ of intervals of size $\delta^{p+2}A^{-r^C}$, $C>7$, such that 
\item{(a)} On $I\in\Lambda_r$, $u^{(r)}(\lambda, m, M)$ is given by a rational function in $(\lambda, m, M)$ of degree at most 
$A^{Cr^3}$.
(Consequently, $\omega^{(r)}$ is rational of degree at most $A^{(2p+1)Cr^3}$ from (N2).)

\item{(b)} For $(\lambda, m, M)\in\bigcup_{I\in\Lambda_r} I$,
$$\aligned\Vert F(u^{(r)})\Vert&<\kappa_r,\\
\Vert \partial F(u^{(r)})\Vert&<\tilde\kappa_r,\endaligned\tag b1$$ 
with $\kappa_{r+1}\ll\kappa_{r}$ and $\tilde\kappa_{r+1}\ll\tilde\kappa_{r}$; and 
$$\aligned \Vert j\cdot\lambda\Vert _{\Bbb T}&\geq \frac{\kappa}{|j|^{\rho}},\, 0<\kappa<1\, \rho>3,\\
\Vert n\cdot\omega^{(r)}\Vert _{\Bbb T}&\geq  \frac{\kappa}{|n|^{\rho}},\, 0<\kappa<1, \,\rho>3,\endaligned\tag b2$$
for $j\in [-A^r, A^r]^2\backslash\{0\}$, and $n\in [-A^r, A^r]^2\backslash\{0\}$.

\item{(c)} Let $N=A^r$. For $(\lambda, m, M)\in\bigcup_{I\in\Lambda_r} I$, $T=T(u^{(r-1)}):=F'(u^{(r-1)})$ satisfies

$\Vert T_N^{-1}\Vert <\delta^{-(p+2)}A^{r^C}$,

$|T_N^{-1}(j, j')|<\delta^{-(p+2)}e^{-|k-k'|^\gamma}$, for $|k-k'|>r^{C}$,

where $T_N$ is $T$ restricted to $[-N, N]^{4}$.
\item{(d)} Each $I\in\Lambda_r$ is contained in an interval $I'\in\Lambda_{r-1}$ and 
$$\text{meas}(\bigcup_{I'\in\Lambda_{r-1}} I'\backslash \bigcup_{I\in\Lambda_{r}}I)<\frac{\kappa}{r^{5}}\, (0<\kappa<1).$$
\smallskip 
The iteration holds with 

$$\delta_r<\delta^pA^{-(\frac{4}{3})^r}, \, \tilde\delta_r<\delta^pA^{-\frac{1}{2}(\frac{4}{3})^r}; \kappa_r<\delta^{2p+1} A^{-(\frac{4}{3})^{r+2}}, \, \tilde\kappa_r<\delta^{2p+1} A^{-\frac{1}{2}(\frac{4}{3})^{r+2}}.\tag W$$
\bigskip
We remark that the approximate solutions $u^{(r)}$ are defined, a priori, on $\Lambda_r$, but using the 
derivative estimates in (Hii) together with (W), as $C^1$ {\it functions} they can be
extended to the entire $(\lambda, m, M)$-space, 
thus verifying (Hi-iii). We call the intervals in $\Lambda_r$, the {\it good intervals}.
\smallskip
\noindent{\it Remark.} Cf. sect.~3.1 \cite{W3}, which provides some intuitions to a similar set of induction hypothesis. 
\bigskip
\subheading{3.2 Scheme of the proof}

We iteratively solve the $P$ and then the $Q$-equations, with $u^{(0)}(-e_1, h_1)=a_1$, $u^{(0)}(-e_1, h_2)=a_2$ and $0$ otherwise; $v^{(0)}(e_1, -h_1)=a_1$, $u^{(0)}(e_1, -h_2)=a_2$ and $0$ otherwise;
$\omega^{(0)}$ as in (2.2).

It consists of two distinct types of steps. The first ones do direct perturbations; while the second use the $(\theta, \phi)$ estimates in the Main Lemma and the covariance
relation (2.17).
The first steps work directly in the $(\lambda, m, M)$ variables and are relatively straightforward, under the condition in the Theorem: $h_1\not\parallel h_2$. 
A key to the second steps is the bi-rational approximations of the diffeomorphisms in (3.1)-(2.13) on the set of intervals $\Lambda_r$ in (Hiv),
which we will elaborate. 
\bigskip
\subheading{3.3 The initial steps}

Let $F'$ be the linearized operator in (2.8). The initial $R$ steps ($R$ to be determined) directly perturb about the diagonal matrix $D$ in (2.9)
to show that $F'_N$ are invertible for $N=A^r$, $1\leq r\leq R$, after excisions in $(\lambda, m, M)$. The matrix elements of $D$ are 
polynomials in $(\lambda, m, M)$. So we only need to exclude the matrix elements being identically zero. The following Lemma serves this purpose.

\proclaim{Lemma 3.1} If 
$$h_1 \not\parallel h_2,$$ 
then the matrix elements in (2.9),
$$D(\pm n, j)\not\equiv 0,\tag 3.3$$ 
for all $(\pm n, j)\in\Bbb Z^4\backslash{\{(\mp e_k, \pm h_k)\}_{k=1}^2}$, i.e., 
on the domain of the $P$-equations.
\endproclaim

\demo{Proof}
Clearly it suffices to take the positive sign and work with
$D(n,j)$. 
Toward that end, write the {\it fixed}
$$0\neq h_1:=(x, y),$$
$$0\neq h_2:=(x', y'),$$
and the {\it generic}
$$(n, j)=(n_1, n_2, j_1, j_2)\in\Bbb Z^4.$$ 
Writing out the polynomial for $D:=D(n, j)$, one has
$$\aligned D=&n_1[(x\lambda_1+y\lambda_2+m)^2+M]+n_2[(x' \lambda_1+y'\lambda_2+m)^2+M]+ (j_1\lambda_1+j_2\lambda_2+m)^2+M\\
=&(j_1^2+n_1 x^2+n_2{x'}^2)\lambda_1^2+(j_2^2+n_1 y^2+n_2{y'}^2)\lambda_2^2+2(j_1j_2+n_1 xy+n_2x'y')\lambda_1\lambda_2\\
&+2m[n_1(x\lambda_1+y\lambda_2)+n_2(x'\lambda_1+y'\lambda_2)+(j_1\lambda_1+j_2\lambda_2)]+(n_1+n_2+1)(m^2+M).\endaligned\tag 3.4$$
Using the last term, clearly when
$$n_1+n_2+1\neq 0,$$
$$D\not\equiv 0.$$
When $$n_1+n_2+1=0,\tag 3.5$$ 
there are the three additional equations, if the first three coefficients are null:
$$\align & j_1^2+n_1 x^2+n_2{x'}^2=0,\tag 3.6\\
&j_2^2+n_1 y^2+n_2{y'}^2=0,\tag 3.7\\
&j_1j_2+n_1 xy+n_2x'y'=0. \tag 3.8\endalign$$
Viewing $n_1, n_2, j_1$ and $j_2$ as the variables, we use the four equations (3.5)-(3.8) to solve for them. 
From (3.5), 
$$n_2=-n_1-1.\tag 3.9$$
Substituting  (3.9) in (3.6)-(3.8), gives 
$$\align &j_1^2={x'}^2-n_1(x^2-{x'}^2),\tag 3.10\\
&j_2^2={y'}^2-n_1(y^2-{y'}^2),\tag 3.11\\
&j_1j_2=x'y'-n_1(xy-x'y').\tag 3.12\endalign$$
Multiplying (3.10) and (3.11) and squaring (3.12) yields
$$[{x'}^2-n_1(x^2-{x'}^2)][y^2-n_1(y^2-{y'}^2)]=[x'y'-n_1(xy-x'y')]^2,\tag 3.13$$
which leads to 
$$n_1^2(xy'-x'y)^2+n_1(xy'-x'y)^2=0.\tag 3.14$$
Since $h_1\not\parallel h_2$, 
$$xy'-x'y\neq 0,$$
we have that the solutions to (3.14), hence (3.5)-(3.8) are
$$n_1=0, n_2=-1, j_1=\pm x', j_2=\pm y',\tag 3.15$$
or
$$n_1=-1, n_2=0, j_1=\pm x, j_2=\pm y.\tag 3.16$$
Finally using the coefficient in front of $m$ in the fourth term in the second equality in (3.4), gives that in
order for {\it all} five coefficients to be null, one needs to take the plus sign solutions in (3.15) and (3.16).
Clearly reversing the signs gives the solutions to $D(-n, j)\equiv 0$ and concludes the proof. 
\hfill $\square$
\enddemo

Initialize $u^{(0)}$ to be 
$$u^{(0)}=a_1\delta_{(-e_1, h_1)}+a_2\delta_{(-e_1, h_2)}\tag 3.17$$
from (2.3),
$$v^{(0)}=a_1\delta_{(e_1, -h_1)}+a_2\delta_{(e_1, -h_2)}\tag 3.18$$
and
$$\omega^{(0)}=((h_1\cdot\lambda+m)^2+M, (h_2\cdot\lambda+m)^2+M)\tag 3.19$$
from (2.2). Define $u^{(r)}$, $v^{(r)}$ and $\omega^{(r)}$ for $r\geq 1$ as in (N1)-(N3). 

The following is a direct consequence of Lemma~3.1.

\proclaim{Lemma 3.2a} If 
$$h_1 \not\parallel h_2,$$ 
(Hi-iii), (Hiv, a, b1, c, d) and (W) are verified for $1\leq r\leq |\log\delta|^{3/4}$. 
\endproclaim

\noindent{\it Remark.} In sect.~3.4, Lemma~3.2 will complete Lemma~3.2a by proving (Hiv, b2) for $1\leq r\leq |\log\delta|^{3/4}$.
\smallskip
\demo{Proof of Lemma~3.2a} Since $h_1\not\parallel h_2$, Lemma 3.1 is available.
Set 
$$|D(\pm n, j)|>2\delta ^{p/2},\tag 3.20$$ 
where $D(\pm n, j)$ as in (3.4). Using that $D(\pm n, j)$ are 
quadratic polynomials in $(\lambda, m, M)$ gives that away from
a set of measure less than $\delta^{p/5}$,
all diagonal elements satisfy (3.20) for $1\leq r\leq |\log\delta|^{3/4}$.
The details are as follows.

We use the first expression for $D$ in (3.4). There are two cases.

$\bullet$ $n_1+n_2+1\neq 0$

Since $\partial D/\partial M=n_1+n_2+1$, 
$$|\partial D/\partial M|\geq 1,$$ leading to measure excision $\ll \delta^{p/3}$.

$\bullet$ $n_1+n_2+1=0$

We make a change of variables and set 
$$\aligned X_1&=x\lambda_1+y\lambda_2+m,\\
X_2&=x'\lambda_1+y'\lambda_2+m,\\
X_3&=j_1\lambda_1+j_2\lambda_2+m.\endaligned$$
The matrix of change of variables:
$$J=\pmatrix x&y&1\\x'&y'&1\\j_1&j_2&1\endpmatrix,$$
$\det J\in\Bbb Z$.

(a) $\det J\neq 0$, so $|\det J|\geq 1$.
$$D=n_1X_1^2+n_2X_2^2+X_3^2$$
and 
$$\frac{\partial^2 D}{\partial X_3^2}=2,$$
leading to measure excision in $(X_1, X_2, X_3)$ $\ll \delta^{p/5}$, hence also 
in $(\lambda, m)$, and finally $(\lambda, m,M)$.

(b) $\det J=0$. Since $h_1\not\parallel h_2$, $X_1$ and $X_2$ are linearly independent,
so $$X_3=c_1X_1+c_2X_2,$$
with $c_1+c_2=1$ (from the $m$-component),
and 
$$\aligned D=&n_1X_1^2+n_2X_2^2+(c_1X_1+c_2X_2)^2\\
=&(n_1+c_1^2)X_1^2+(n_2+c_2^2)X_2^2+2c_1c_2X_1X_2.\endaligned$$
Now $c_1$ and $c_2$ must satisfy 
$$\aligned j_1&=c_1x+c_2x'\\
 j_2&=c_1y+c_2y'.\endaligned$$
 So 
 $$c_1=\frac{\det\pmatrix j_1&x'\\j_2&y'\endpmatrix}{\det\pmatrix x&x'\\y&y'\endpmatrix},$$
 and 
  $$c_2=\frac{\det\pmatrix x&j_1\\y&j_2\endpmatrix}{\det\pmatrix x&x'\\y&y'\endpmatrix}.$$
  Since $$\det\pmatrix x&x'\\y&y'\endpmatrix=\ell\neq 0,$$
  and all the determinants are integer valued,
  $c_1$ and $c_2\in\Bbb Z/|\ell|$. 
  
  If both $n_1+c_1^2=0$ and $n_2+c_2^2=0$, then since $n_1+n_2+1=0$,
  this yields 
  $$c_1^2+c_2^2=1.$$
  Since $$c_1+c_2=1,$$
  this means that $c_1=1$ and $c_2=0$, leading to
  $(n_1, n_2)=(-1, 0)$ and $(j_1, j_2)=(x, y)=h_1$,
  which is in the domain of the $Q$-equations, therefore
  excluded; or $c_1=0$ and $c_2=1$, leading to
   $(n_1, n_2)=(0, -1)$ and $(j_1, j_2)=h_2$,
   which are again excluded. 
  
  Otherwise, either $n_1+c_1^2\neq 0$ or $n_1+c_2^2\neq 0$. 
  Assume $n_1+c_1^2\neq 0$. Since 
  $$n_1+c_1^2\in\Bbb Z/|\ell|^2,$$
  $$|n_1+c_1^2|\geq \frac{1}{|\ell|^2}=\Cal O(1).$$
  So 
  $$\frac{\partial^2 D}{\partial X_1^2}\geq \frac{2}{|\ell|^2},$$
  yielding measure excision $\ll \delta^{p/5}$.
  Clearly the same holds in the sub-index $2$ variables and concludes the 
  measure estimates.
  \smallskip
 Direct perturbation theory using (3.20), small $\delta$, and 
$$r\leq |\log\delta|^{3/4}:=R,\tag 3.21$$ 
then leads to (Hi-iv) and (W), except the Diophantine estimates (Hiv, b2), which are not needed for the initial scales.
(Cf. the proofs of Lemma~3.1 and the Corollary in \cite{W3}.)

Let us emphasize nonetheless, two points in the proof:

\item{(i)}  The interval size is determined by the size of perturbations which preserves the estimates in (Hiv,~c). For the 
first $R$ steps, it simplifies to perturbations of the diagonals, which give the size to be of order (at least) $\delta^{p+2}$
from (3.20). The inverse of the size of the intervals gives clearly,  an upper bound on the number of intervals at each stage $r$.
(This upper bound plays an important role in sect.~3.3.)

\item{(ii)} From (2.11), $u^{(r)}$ is rational, and 
$$\text{deg }u^{(r+1)}\lesssim [A^{(r+1)}]^4\text {deg } u^{(r)},$$
where the volume factor comes from the inverse of the linearized operator. This leads to $$\text{deg }u^{(r)}\leq A^{Cr^3}.\tag 3.22$$

\hfill$\square$
\enddemo

Clearly the arguments in (i) and (ii) hold in greater generality, and will be used for $r\geq R$ as well.   

\bigskip
\subheading{3.4 Bi-rational approximations and the general steps}

Increasing from scale $A^r$ to $A^{r+1}$ for $r\geq R$, the general idea is 
to pave $[-A^{r+1}, A^{r+1}]^4$ by $[-A^{r}, A^{r}]^4$ and (much) smaller cubes of size $\sim r^C$, $C>1$ (cf., Lemma 3.4), at distances at least $A^{r}/2$ from the origin. 

There are two types of $r^C$-cubes. Type 1 uses the Main Lemma and then the projection Lemma~2.2; Type 2 makes direct excisions using the Lipschitz families of functions constructed in Lemma~B, 
satisfying (2.24) and (2.25). Type 1 is more subtle.  

Recall that the Main Lemma is applied {\it per interval}.
On Type 1 $r^C$-cubes, it is therefore essential to keep the complexity, i.e., number of good intervals, low. 
Using (W), this could be accomplished by replacing $u^{(r)}$ by appropriate $u^{(r')}$, with $r'\ll r$.
(In fact this is why we have chosen these much smaller cubes.)
Subsequently, the main work is to verify that on $I\in\Lambda_{r'}$,
one may indeed approximate the diffeomorphism defined by (3.1)-(N2), by bi-rational maps. 

\noindent{\it Bi-rational approximations.}

The map $\Cal M_r$: $(\lambda, m, M)\mapsto (\lambda, \omega)$ is clearly rational on $I\in\Lambda_r$ using (Hiv,~a). 
Below we show that $\Cal M_r^{-1}$ admits rational approximations $\mu_r$ on $\Cal M_r(I)$.  

\proclaim{Lemma 3.3} Let $I\in\Lambda_r$, $r\in \{0, 1, 2, ...\}$. There exist rational maps $\mu_r$ on $\Cal M_r(I)$, such that 
$$\mu_r\Cal M_r=\Bbb I+\Cal O(\delta_r^2),\tag 3.23$$
and 
$$\Cal M_r\mu_r=\Bbb I+\Cal O(\delta_r^2),\tag 3.24$$
for $r\in \{0, 1, 2, ...\}$.
\endproclaim

\demo{Proof} At the $r^{\text {th}}$ step, the $Q$-equation reads:
$$-\omega_k+(h_k\cdot\lambda+m)^2+M+\delta^{2p} \frac{[{(u^{(r)}*v^{(r)}})^{*p}*u^{(r)}] (\lambda, m, M)} {a_k}(-e_k, h_k)=0,\, k=1, 2,\tag 3.25$$
and defines the map 
$$\Cal M_r: (\lambda, m, M)\mapsto (\lambda, \omega),$$
on $(0, 2\pi)^4$. 

Define 
 $$\Gamma=\{\lambda|(h_1-h_2)\cdot\lambda>\kappa',\, 0<\kappa'<1\}.$$ 
Using (Hii), the matrix of change of variables
$$(\lambda, m, M)\mapsto (\lambda, \omega)$$
 is
$$J=\pmatrix 1&0&2(h_1\cdot\lambda+m)h_{1,1}&2(h_2\cdot\lambda+m)h_{2,1}\\
 0&1&2(h_1\cdot\lambda+m)h_{1,2}&2(h_2\cdot\lambda+m)h_{2,2}\\
 0&0&2(h_1\cdot\lambda+m)+\Cal O(\delta^{2p})&2(h_2\cdot\lambda+m)+\Cal O(\delta^{2p})\\
0&0&1&1\endpmatrix+\Cal O(\delta^{2p}).$$
The Jacobians satisfy, 
$$\aligned \det J&=2(h_1-h_2)\cdot\lambda+\Cal O(\delta^{2p})=\Cal O(1),\\
 \det J^{-1}&=\frac{1}{2(h_1-h_2)\cdot\lambda}\big[1+\frac{\Cal O(\delta^{2p})}{2(h_1-h_2)\cdot\lambda}\big]=\Cal O(1).\endaligned \tag 3.26$$
 on $\Gamma\times (0, 2\pi)^3$.
Hadamard's global inverse function theorem, see sect.~6.2 \cite {KP}, gives that
the inverse
$$\Cal M^{-1}_r: (\lambda, \omega)\mapsto (\lambda, m, M),$$
is well defined on  $\Gamma\times (0, 2\pi)^2$. 



Denote $I\in\Lambda_r$ by $I_r$. Set $I_0=\Gamma\times (0, 2\pi)^2$. 
From (Hiv, d), there exist 
$$I_{r-1}\in\Lambda_{r-1},\, I_{r-2}\in\Lambda_{r-2}, ...,I_0\subset (0, 2\pi)^4,$$ 
with the nested property: 
$$I_r\subset I_{r-1}\subset I_{r-2}\subset \cdots \subset I_0\subset  (0, 2\pi)^4,$$
on which (Hiv) hold at steps $i$, $i\leq r$ respectively. It follows that 
$$|I_i|\sim O(A^{-i^C})\gg \delta_i,$$
using (W), for all $i\leq r$. 

Let $$\Cal I_i=\Cal M_i(I_i),\, i\leq r.$$
From (3.26)
$$|\Cal I_i|\sim A^{-i^C}\gg \delta_i.$$
Using (3.1) and (3.2), one may therefore assume
$$\Cal I_r\subset \Cal I_{r-1}\subset \Cal I_{r-2} \subset \cdots\subset\Cal I_0:=\Cal M_0 (I_0),\tag 3.27$$
on which (Hiv) hold in the $(\lambda, \omega)$ variables.
 
Below we make rational approximation to $\Cal M_r^{-1}$ on $\Cal I_r$.
When $r=0$, the solutions to (3.25) are given by (2.35)-(2.36):   
$$\aligned m=m^{(0)}(\lambda, \omega) &= \frac{\Omega_1-\Omega_2-(h_1\cdot\lambda)^2+(h_2\cdot\lambda)^2}{2(h_1-h_2)\cdot\lambda},\\
M=M^{(0)}(\lambda, \omega)&=\frac{\Omega_1+\Omega_2}{2}-\frac{[(h_1-h_2)\cdot\lambda]^2}{4}-\frac{(\Omega_1-\Omega_2)^2}{4[(h_1-h_2)\cdot\lambda]^2},\endaligned\tag 3.28$$
where $$\Omega_k=\omega_k-\delta^{2p}\frac{{(u^{(0)}*v^{(0)})^{*p}*u^{(0)}}}{a_k}(-e_k, h_k),\, k=1, 2,\tag 3.29$$
provided $$(h_1-h_2)\cdot\lambda\neq 0.\tag 3.30$$
Clearly $m^{(0)}$ and $M^{(0)}$ are {\it rational} in $(\lambda, \omega)$. 

Denote the left side of (3.25) by $\Cal Q_r$. When $r=0$, exceptionally, $\Cal Q_0(m^{(0)}, M^{(0)})=0$.
When $r=1$,
$$|\Cal Q_1(m^{(0)}, M^{(0)})|\lesssim \Vert\Delta u^{(1)}\Vert\lesssim\delta_1.$$
The linearized operator $\Cal Q'_1$ is a $2\times 2$ matrix: 
$$\Cal Q'_1=\pmatrix \frac{\partial \Cal Q_{1,1}}{\partial m} &\frac {\partial \Cal Q_{1,1}}{\partial M} \\ \frac{\partial \Cal Q_{1,2}}{\partial m} & \frac{\partial \Cal Q_{1,2}}{\partial M} \endpmatrix
=\pmatrix 2(h_1\cdot\lambda+m)+\Cal O(\delta^{2p})&1+\Cal O(\delta^{2p})\\2(h_2\cdot\lambda+m)+\Cal O(\delta^{2p})&1+\Cal O(\delta^{2p})\endpmatrix.$$
$$\det \Cal Q'_1=2(h_1-h_2)\cdot\lambda+\Cal O(\delta^{2p})=\Cal O(1),$$
from (3.25). So 
$$\Vert [\Cal Q'_1]^{-1}\Vert =\Cal O(1).$$
(Clearly, using (Hii) the above estimate holds for all $[\Cal Q'_i]^{-1}$, $i\leq r$.)
Set 
$$(\Delta m^{(1)}, \Delta M^{(1)}):=[\Cal Q'_1]^{-1}(m^{(0)}, M^{(0)})\Cal Q_1(m^{(0)}, M^{(0)}),$$
and $$(m^{(1)}, M^{(1)}):=(m^{(0)}, M^{(0)})+(\Delta m^{(1)}, \Delta M^{(1)}).$$
The approximate solutions $m^{(1)}$ and $M^{(1)}$ are rational on $\Cal I_1$ and
 $$\Cal Q_1(m^{(1)}, M^{(1)})=\Cal O(\delta_1^2).$$
 (For simplicity, we have written $\Cal O$ for $\Cal O_\delta$, as  $\delta$ is now fixed.)
 
 Assume at step $(i-1)$, $i\geq 1$, there are rational $m^{(i-1)}$ and $M^{(i-1)}$ on $\Cal I_{i-1}$ such that
 $$\Cal Q_{i-1}(m^{(i-1)}, M^{(i-1)})=\Cal O(\delta_{i-1}^2).\tag 3.31$$
 Then 
 $$\Cal Q_{i}(m^{(i-1)}, M^{(i-1)})=\Cal O(\delta_{i})+\Cal O(\delta_{i-1}^2)\sim (1+\delta_i^{1/2}) \Cal O(\delta_{i}).$$
 Set 
$$(\Delta m^{(i)}, \Delta M^{(i)}):=[\Cal Q'_{i}]^{-1}(m^{(i-1)}, M^{(i-1)})\Cal Q_{i}(m^{(i-1)}, M^{(i-1)}),$$
and $$(m^{(i)}, M^{(i)}):=(m^{(i-1)}, M^{(i-1)})+(\Delta m^{(i)}, \Delta M^{(i)}).$$
Then $m^{(i)}$ and $M^{(i)}$ are rational on $\Cal I_{i}$. (Here we used $\Cal I_{i}\subset \Cal I_{i-1}$.)
Since 
$$\Vert [\Cal Q'_{i}]^{-1}\Vert =\Cal O(1),$$
we obtain 
 $$\Cal Q_{i}(m^{(i)}, M^{(i)})\sim (1+\delta_i^{1/2})^2\delta_{i}^2\sim (1+2\delta_i^{1/2})\delta_{i}^2< (1+\delta_i^{1/3})\delta_{i}^2, \, 1\leq 1\leq r.\tag 3.32$$
 Iteration using (3.31) and (3.32) from $i=1$ to $r$, yields 
 $$\Cal Q_{r}(m^{(r)}, M^{(r)})<(1+\delta_r^{1/3})\delta_{r}^2<2\delta_r^2.$$
 
 It follows from the construction and (3.22) that 
 $$\text{deg }m^{(r)},\, \text{deg }M^{(r)} \sim \text {deg } u^{(r)}\sim A^{(2p+1)Cr^3},\tag 3.33$$
 in $(\lambda, \omega)$. 
Defining 
 $$\mu_r=(m^{(r)}, M^{(r)}),$$
concludes the proof.
\hfill$\square$
\enddemo
\bigskip
With Lemma~3.3 in hand, we are ready to make the general steps.
So let us complete the initial steps now.
\proclaim{Lemma 3.2} 
The Diophantine estimates (Hiv, b2) are verified for $1\leq r\leq |\log\delta|^{3/4}$. 
Hence, (Hi-iv) are verified for $1\leq r\leq |\log\delta|^{3/4}$, if $h_1\not\parallel h_2$. 
\endproclaim

\demo{Proof} Using (3.1), (3.2) and (3.26), (Hiv, b2) are satisfied, after excising a set of measure at most
$\kappa$. (Note that diffeormorphism suffices for the measure estimates, and Lemma~3.3 is not used here.)
Taking into account Lemma~3.2a, concludes the proof.
\hfill$\square$
\enddemo

Assume (Hi-iv) hold at stage $r$. To construct $u^{(r+1)}$, we need to control 
$$T_N^{-1}(u^{(r)}) \text { with }N=A^{r+1}.$$
This requires another excision in $(\lambda, m, M)$, which will lead to the next set of intervals $\Lambda_{r+1}$. 

To simplify notations, given two sets of intervals $Z_1$ and $Z_2$, we say that 
$$Z_2\subset Z_1,$$
if for all $I\in Z_2$, there exists $I'\in Z_1$, such that
$I\subset I'$.  We also define
$$\text{meas } (Z_1\backslash Z_2)=\text {meas }(\bigcup_{I'\in Z_1} I'\backslash\bigcup_{I\in Z_2} I.)$$
\bigskip
\subheading{3.5 The new set of intervals for the $(\theta, \phi)$ estimates}
\smallskip
Assume that (Hi-iv) hold at step $r$. To prove that they hold at step $(r+1)$, we
cover $[-A^{r+1}, A^{r+1}]^{4}$ by $[-A^{r}, A^{r}]^{4}$ and smaller
cubes $[-A_0, A_0]^{4}+J$, with $A^r/2<|J|<A^{r+1}$ and $A_0\ll A^{r+1}$.
The following Lemma provides $(\theta,\phi)$ estimates on the $A_0$-cube centered at the origin. 

\proclaim{Lemma 3.4} At stage $r$, scale $N=A^r$, set 
$$A_0=(\log N)^C=r^C(\log A)^C,\tag 3.34$$ with $C>7/c$, and $c$ as in (2.44);
define $$r_0=\frac{\log A_0}{\log A}$$ and 
$$\tilde r_0:=r_0\frac{\log A}{\log 4/3}\simeq C\frac{\log r}{\log 4/3}\ll r.\tag 3.35$$
Then there is  $\Lambda'_{r+1}\subset \Lambda_r$, satisfying
$$\text{meas } (\Lambda_r\backslash\Lambda'_{r+1})<\frac{\kappa}{ r^6},$$
so that on the intervals in the set  $\Lambda'_{r+1}$, (Hiv, b2) hold and there are the following estimates:
$$\aligned \Vert T_{A_0}^{-1}(\lambda, m, M; u^{(\tilde r_0)},v^{(\tilde r_0)} )(\theta,\phi)\Vert &<e^{A_0^\sigma},\\
|T_{A_0}^{-1}(\lambda, m, M; u^{(\tilde r_0)},v^{(\tilde r_0)} )(\theta,\phi)(k, k')| &< e^{- |k-k'|^\gamma},\endaligned\tag 3.36$$
for all $k, k'$ such that $|k-k'|>A_0/10$, provided $(\theta, \phi)$ is in the complement of a set $\Theta_{A_0}(\theta, \phi)\in \Bbb R^2$,
whose sectional measures satisfy
$$\text{meas } [\theta| \text {fixed }\phi, (\theta, \phi)\in \Theta_{A_0}]\leq e^{-A_0^\tau}\, (\tau>0),$$
and $$\text{meas } [\phi|\text {fixed }\theta, (\theta, \phi)\in \Theta_{A_0}]\leq e^{-A_0^\tau} \, (\tau>0).\tag 3.37$$
\endproclaim

\demo{Proof} 
Assume (Hiv, b2) hold at stage $r$, then (Hiv, b2) hold at stage $r+1$ 
after excision of a set of measure at most $A^{-r}$.

For $r$ such that $r_0\leq R$, perturbing directly about the diagonals $D_{\pm}$ defined in (2.16)
leads to (3.36)-(3.37). There is no additional excision in $(\lambda, m, M)$. 

For $r_0>R$, the strategy of the proof is similar to that of Lemma~3.2 in \cite{W3}, and uses
the Main Lemma. (Cf. also sect.~2.5 \cite{W3}.)
Set $$N_1=(\log A_0)^C\simeq(\log r)^C,$$ in view of  (3.34), 
$$\tilde r=\frac{\log N_1}{\log A}\simeq \log\log r\ll r_0\ll r,$$
and 
$$\tilde{\tilde r}=\tilde r \frac{\log A}{\log 4/3}\simeq \log\log r\ll \tilde r_0\ll r.$$

At scale $N_1$, using (Hii, iii), it suffices to replace $u^{(\tilde r_0)}$ by $u^{(\tilde{\tilde r})}$.
To apply the Main Lemma, fix $I\in \Lambda_{\tilde{\tilde r}}$.
By the choice of $\tilde{\tilde r}$, using (W) on $\Lambda_{\tilde{\tilde r}}\cap \Lambda_{\tilde r_0}$ and Lemma ~3.3, we have
$$\aligned \Vert T_{N_1}(\lambda, m^{(\tilde{\tilde r})}, M^{(\tilde{\tilde r})}; u^{(\tilde{\tilde r})}, v^{(\tilde{\tilde r})})(\theta,\phi)-T_{N_1}(\lambda, m, M; u^{(\tilde r_0)}, v^{(\tilde r_0)})(\theta, \phi)\Vert 
&\leq \Cal O(\delta_{\tilde{\tilde r}})+\Cal O(\delta_{\tilde{\tilde r}}^2)\\
&\leq e^{-\tilde\gamma N_1},\endaligned\tag 3.38$$
for some $\tilde\gamma>0$, and $|\theta|\leq 2e^{2(\log A_0)^2}=2e^{N_1^{2/C}}$, $|\phi|\leq 2e^{(\log A_0)^2}=2e^{N_1^{1/C}}$. 

(This follows the restrictions on $\theta$ and $\phi$ in region (i),
with $A_0$ in place of $N$, in the Main Lemma. 
Recall that otherwise, Diophantine conditions (Hiv, b2) suffice and there is no additional excisions in $(\lambda, \omega)$ or $(\lambda, m, M)$.
Here ``$\cap$" is in the sense of intersections of the intervals in the two sets; note from (Hiv,~d) that each interval in $\Lambda_{\tilde r_0}$ is {\it contained} in an interval 
in $\Lambda_{\tilde{\tilde r}}$.)

The set on which (3.36) holds at scale $A_0$ are obtained after excisions on the set at scale $N_1$ 
leading to a set $\Omega'_3$, the image set of $\Omega_3$ in (2.52) under the map $\Cal M_{\tilde r_0}^{-1}$.

Rename $\Omega'_3$, $\Omega'_{3, I}$. Set 
$$\bar \Gamma_{A_0}=\cap_I \Omega'_{3, I},$$
and 
$$\Gamma_{A_0}=\bar \Gamma_{A_0}\cap\Lambda_{\tilde r_0}.$$
Define 
$$\Lambda'_{r+1}=\Lambda_{r}\cap \Gamma_{A_0}.$$


There are at most $$A^{(\tilde{\tilde r})^C}\simeq A^{(\log\log r)^C}$$
such intervals $I$ in $\Lambda_{\tilde{\tilde r}}$, by using (Hiv) at stage 
$$\tilde{\tilde r}\sim\log\log r\ll r.$$ 
So using (2.52) and (3.26), and adding also the Diophantine excisions, we have 
$$\aligned \text{meas } (\Lambda_r\backslash\Lambda'_{r+1})<A^{(\log\log r)^C}A_0^{-c}+A^{-r}&<A^{(\log\log r)^C}r^{-Cc}\\
&<\frac{\kappa}{r^6},\endaligned\tag 3.39$$
if $Cc>7$.
On $\Lambda'_{r+1}$, (Hiv, b2), (3.36) and (3.37) hold.
\hfill$\square$
\enddemo
\bigskip
\subheading{3.6 The new set of good intervals}

We first convert the $(\theta,\phi)$ estimates in (3.36) for the $A_0$-cube centered at the origin to 
$A_0$-cubes centered at large $J\in\Bbb Z^{4}$ at $(\theta,\phi)=0$. 
We use the projection Lemma~2.2 in sect.~2.7, as well as the families of Lipschitz functions in Lemma~B, sect.~2.4,
for this purpose.
\smallskip
\proclaim {Lemma~3.5} 
There exists $\Lambda_{r+1}\subset\Lambda'_{r+1}\subset\Lambda_r$, 
satisfying $$\text{meas } (\Lambda_r\backslash \Lambda_{r+1})<\frac{\kappa}{r^5},\tag 3.40$$
provided $C>max (1/\tau, 7/c)$. 
On the intervals in the set $\Lambda_{r+1}$, (Hiv, b2) hold and $T_{[-A_0, A_0]^{2d}+J}^{-1}(u^{(\tilde r_0)})$,
with $A_0$ and $\tilde r_0$ as in (3.34) and (3.35) respectively,
satisfy the upper bounds in (3.36), for all $J$ with $A^r/2<|J|<A^{r+1}$.
\endproclaim

\demo{Proof} Fix $I\in\Lambda'_{r+1}$. From (3.36) of Lemma~3.4, on the cube $[-A_0, A_0]^4$, there are the estimates:
$$\aligned \Vert T_{A_0}^{-1}(\lambda, m, M; u^{(\tilde r_0)},v^{(\tilde r_0)} )(\theta,\phi)\Vert &\lesssim e^{A_0^\sigma},\\
|T_{A_0}^{-1}(\lambda, m, M; u^{(\tilde r_0)},v^{(\tilde r_0)} )(\theta,\phi)(k, k')| &\lesssim e^{-|k-k'|^\gamma},\endaligned \tag 3.41$$
for all $k, k'$ such that $|k-k'|>A_0/10$, away from a set in $(\theta, \phi)$ of sectional measures at most $e^{-A_0^\tau}$.

Divide $(\theta, \phi)$ into 
two regions, as in sect.~2: 
\item {(i)} $|\theta|\leq 2e^{2(\log A_0)^2},\,\, |\phi|\leq 2e^{(\log A_0)^2}$;
\item {(ii)} otherwise.
     
Write $$J= (N_1, N_2; J_1, J_2)\in\Bbb Z^4, \, A^r/2\leq |J|\leq A^{r+1}.$$
We work on $\Cal M_{\tilde r_0}(I):=\Cal I$ in the $(\lambda, \omega)$ variables, where $\Cal M_{\tilde r_0}$ as defined in (2.13)-(3.1).
Set $$\aligned &\theta=N_1\omega_1+N_2\omega_2,\\
&\phi=J_1\lambda_1+J_2\lambda_2.\endaligned$$
Clearly for any given $J$, $(\theta, \phi)$ as defined above falls into one of the two regions.
\bigskip
We start with region (ii). On region (ii),
at least one of the following holds:
$$ |\theta|>2e^{2(\log A_0)^2},$$
or $$|\phi|>2e^{(\log A_0)^2}. $$
Further one may assume 
$$|\,|\phi|-\sqrt{|\theta|}\,|<e^{(\log A_0)^2},\tag 3.42$$
as other wise $T_\Lambda$ is invertible.

If $ |\theta|>2e^{2(\log A_0)^2}$, (3.42) leads to 
$$ \aligned&|\theta|>2e^{2(\log A_0)^2}\\
&|\phi|<\sqrt{|\theta|}+e^{(\log A_0)^2}.\endaligned\tag ii'$$
Similarly, if $|\phi|>2e^{(\log A_0)^2}$, (3.42) yields
$$\sqrt{|\theta|}>|\phi|-e^{(\log A_0)^2},$$
hence
$$|\theta|>e^{2(\log A_0)^2}.$$
This leads to 
$$ \aligned&|\theta|>e^{2(\log A_0)^2}\\
&|\phi|<\sqrt{|\theta|}+e^{(\log A_0)^2}.\endaligned\tag ii''$$
Therefore in region (ii), one may assume that 
for $$\Lambda=[-A_0, A_0]^4+J, \, A^r/2\leq |J|\leq A^{r+1},$$
$$\Cal O(1) A^{r+1}>|\theta|=|N_1\omega_1+N_2\omega_2|>e^{2(\log A_0)^2}.\tag 3.43$$
\bigskip
Fix $\lambda=(\lambda_1, \lambda_2)$.
We use the set of Lipschitz functions $\theta_i$ in Lemma~B, to make direct excisions.
Restricting to $\omega$ given by (2.13)-(3.1), (2.24) yields
$$\aligned \Vert \theta_i\Vert_{\text{Lip}(\omega)}\Vert &\leq \Vert \theta_i\Vert_{\text{Lip}(\omega)}+\Cal O(1)\Vert \theta_i\Vert_{\text{Lip}(M)}+\Cal O(1)\Vert \theta_i\Vert_{\text{Lip}(m)}\\
&\leq\Cal O(A_0|\phi|),\endaligned\tag 3.44$$
where $\Vert \, \Vert_{\text{lip}}$ denotes the restricted Lipschitz norm, and we used (3.26).  

To verify (3.41), it suffices to require
$$|\theta-\theta_i|>2e^{-A_0^\sigma}, \tag 3.45$$
for all $i$, $1\leq i\leq 2(2A_0+1)^4$.

Substitute (3.43) for the $\theta$ in (3.45); use (ii') and (ii'') in (3.44). Varying $\omega$ and
then using (3.26) again, (3.45) 
excises a set in $(\lambda, m, M)$ of measure at most 
$$\Cal O(1)e^{-A_0^\sigma} (A^{r+1})^4(2A_0+1)^4\leq  e^{-A_0^\tau}.\tag 3.46$$
In the complement of the set (3.41) hold.
(Cf. the related proof of Lemma~3.4 in \cite{W1}.)
\bigskip
In region (i), make a partition, so that  
$$\aligned [-2e^{2(\log A_0)^2}, 2e^{2(\log A_0)^2}]\times [-2e^{(\log A_0)^2}, 2e^{(\log A_0)^2}]&=\cup_K\{[-1, 1]^2+K\}\\
&:=\cup_K L_K,\endaligned$$
where $K\in\Bbb Z^2$, satisfying  $$0\leq |K|\leq  10e^{2(\log A_0)^2}.\tag 3.47$$
Fix a $K$ and let  
$$\Cal I'\subset \Cal I,$$ 
be a two dimensional section of $\Cal I$, after
fixing two of the coordinates in $(\lambda_1, \lambda_2, \omega_1, \omega_2)$. 

Set 
$$\Cal S=\Cal I'\times \{\Theta_{A_0}\cap L_K\}\subset \Bbb R^{4}.$$
We have
$$\text{meas }\Cal S\lesssim e^{-A_0^\tau},$$
for all such $\Cal I'$.

The set $\Cal S$ is described by the opposite of (3.36). 
Replacing the $\ell^2$ norm by the Hilbert-Schmidt norm, and $\Cal M^{-1}_{\tilde r_0}$ by its rational approximation 
$\mu_{\tilde r_0}$ in Lemma 3.3, satisfying (3.23) and (3.24). Since the matrix elements of the inverse 
is the division of two determinants, (3.36) can be expressed as algebraic inequalities
in the matrix elements of degree at most $A_0^C$. Using that each matrix element is linear in $\theta$, quadratic in $\phi$ and at most of degree $e^{(C{\log A_0})^3}$ in $\lambda$ and $\omega$,
$\Cal S$ is of degree at most 
$$\text{deg }\Cal S\leq A_0^Ce^{(C{\log A_0})^3}\lesssim e^{({\log A_0})^4}.$$ 
So
$$\log\text{deg }\Cal S\ll \log|J|\ll-\log\text{meas } \Cal S.\tag 3.48$$

Let ``maxv" denote maximum in absolute value. 
Set $$\theta=\text{maxv }( {N_1\omega_1, N_2\omega_2}),$$
and $$\phi=\text{maxv } ({J_1\lambda_1, J_2\lambda_2}).$$
Since 
$$A^r\leq |J|\leq A^{r+1},$$
$$\max (|N_1|, |N_2|; |J_1|, |J_2|)>A^r/2,$$
$$\max (|\theta|, |\phi|)>A^r/2.\tag 3.49$$

From (3.48) and (3.49), we may invoke Proposition~5.1 in \cite{BGS} (and its proof). 
Summing over $K$ satisfying (3.47), we have consequently,  
$$\aligned \Vert T_{\Lambda}^{-1}(\lambda, \omega; u^{(\tilde r_0)},v^{(\tilde r_0)} )\Vert &\lesssim e^{A_0^\sigma},\\
|T_{\Lambda}^{-1}(\lambda,\omega; u^{(\tilde r_0)},v^{(\tilde r_0)} )(k, k')| &\lesssim e^{- |k-k'|^\gamma},\endaligned $$
for all $k, k'$ such that $|k-k'|>A_0/10$, 
after excision of a set in $(\lambda, \omega)$ of measure at most $A^{-2c'r}$, $c'>0$.
Therefore using Lemma ~3.3 and (3.26)
$$\aligned \Vert T_{\Lambda}^{-1}(\lambda, m, M; u^{(\tilde r_0)},v^{(\tilde r_0)} )\Vert &\lesssim e^{A_0^\sigma},\\
|T_{\Lambda}^{-1}(\lambda, m, M; u^{(\tilde r_0)},v^{(\tilde r_0)} )(k, k')| &\lesssim e^{- |k-k'|^\gamma},\endaligned \tag 3.50$$
for all $k, k'$ such that $|k-k'|>A_0/10$, 
after excision of a set in $(\lambda, m, M)$ of measure at most $\Cal O(1) A^{-2c'r}$, $c'>0$.
\bigskip
Adding the measure excised by (3.45) in region (ii), satisfying (3.46), we obtain that (3.50) 
hold for all $$\Lambda=[-A_0, A_0]^4+J, \, A^r/2\leq |J|\leq A^{r+1},$$
after excision of a set $\Cal L'$ in $(\lambda, m, M)$ of measure at most $A^{-c'r}$, $c'>0$.
 
The number of intervals at stage $\tilde r_0$ is bounded above by
$$A^{{\tilde r_0}^C}\simeq A^{(\log A_0)^C}\simeq A^{(\log r)^C},\tag 3.51$$
by using (Hiv) at stage $\tilde r_0$.
Multiplying the measure estimates on $\Cal L'$ by (3.51), 
we obtain that the upper bounds in 
(3.36) hold for all $$\Lambda=[-A_0, A_0]^4+J, \, A^r/2\leq |J|\leq A^{r+1},\tag 3.52$$
after excision of a set $\Cal L$ in $(\lambda, m, M)$ of measure at most 
$$\Cal O(1) A^{(\log r)^C}\cdot A^{-c'r}<A^{-c'r/2},\, c'>0.$$  

Setting $$\Lambda_{r+1}=\Lambda'_{r+1}\backslash \Cal L.$$
concludes the proof.
\hfill$\square$
\enddemo

\noindent{\it Remark.} One may also directly use Lemma ~2.2 to arrive at (3.50),
and view Proposition~5.1 \cite{BGS} as a corollary.
\smallskip
\proclaim {Lemma~3.6} Assume (Hi-iv) hold at stage $r$, on the set $\Lambda_{r+1}$ in Lemma~3.5,
 (Hiv,~c,~d) hold at stage $r+1$.
\endproclaim

\demo{Proof} 
From (3.34) and (3.35), 
$$A_0=r^C(\log A)^C$$ 
and 
$$\tilde r_0=\frac{\log A_0}{\log 4/3}\sim \frac{\log r^C}{\log 4/3}.$$
Using (Hii) between $\tilde r_0$ and $r$, and (W) yields 
$$\Vert u^{(r)}-u^{(\tilde r_0)}\Vert <2\delta_{\tilde r_0}<2e^{-(\log A) A_0}.$$
So the upper bounds in (3.50) are essentially preserved for all  
$T_\Lambda^{-1}(u^{(r)})$, with $\Lambda$ as in (3.52),
where we also used (3.2), $|J|\leq A^{r+1}$ and (Hiii).

Similarly, using (Hii) between $(r-1)$ and $r$, (W), (3.2), $|J|\leq A^{r}$ and (Hiii), give that 
$T_{A^{r}}^{-1}(u^{(r)})$ satisfies essentially the same bound as in
(Hiv, c).

Consequently, resolvent series using $T_{A^{r}}^{-1}(u^{(r)})$ and 
$T_\Lambda^{-1}(u^{(r)})$, proves that (Hiv, c) holds. (This is the same as in (3.18)-(3.25) leading to Lemma~3.5 \cite{W3}.)
Since the estimates on $T_{A^{r}}^{-1}(u^{(r)})$ and $T_\Lambda^{-1}(u^{(r)})$ are stable under perturbations of size $A^{-(r+1)^C}$, this produces
the next set of intervals $\Lambda_{r+1}$ of size $A^{-(r+1)^C}.$ 
(Hiv, d) follows by using the same consideration leading to the Corollary to the Main Lemma, at the end of sect.~2.8.
\hfill $\square$
\enddemo
\bigskip
\subheading{3.7 Proof of the Theorem}
\smallskip
\noindent $\bullet$ {\it Construction of} $u^{(r+1)}$
 
Using Lemma~3.6, we obtain $u^{(r+1)}$ by setting
$$\Delta u^{(r+1)}=-T_{A^{r+1}}^{-1}(u^{(r)})F(u^{(r)}),\tag 3.53$$
and $u^{(r+1)}=u^{(r)}+\Delta u^{(r+1)}$. 
This verifies (Hiv, a) at stage $r+1$. The estimates (Hiv, b1) at stage $r+1$ follow
by direct computation. ((Hiv, b2) hold on $\Lambda_{r+1}$ from Lemma~3.5.)
A standard extension argument 
as in sect.~10, (10.33-10.37)  \cite{B2}, 
shows that (Hi-iii) hold at stage $r+1$ as well. This concludes the induction.
(Cf. the proof of the Corollary in sect.~3.2 \cite{W3}.)
\hfill$\square$
\bigskip

\demo{Proof of the Theorem}
The induction process above solves iteratively the $Q$ and the $P$-equations, with the convergence estimates in (W). The measure estimate (1.2)
follows from (Hiv, d).
The estimates in (1.3) and (1.4) follow from (Hii), (3.1), (3.2) and (W),
and conclude the proof. \hfill$\square$
\enddemo

\enddocument